\documentclass{elsarticle}
\usepackage{color, amsmath, latexsym,amsfonts,amssymb,bm, dsfont}
\usepackage{booktabs}
\usepackage{comment}
\usepackage{graphicx}
\usepackage{subfigure}
\usepackage{float}
\allowdisplaybreaks[4]

\renewcommand{\footnote}{}

\reversemarginpar

\newtheorem{theorem}{Theorem}[section]
\numberwithin{equation}{section}
\newtheorem{expl}{Example}[section]
\newenvironment{example}{\begin{expl}\rm}{\end{expl}}
\newtheorem{prop}[theorem]{Proposition}
\newtheorem{lemma}[theorem]{Lemma}
\newtheorem{corollary}[theorem]{Corollary}

\newtheorem{definition}{Definition}
\newtheorem{rem}[theorem]{Remark}
\newtheorem{assumption}{Assumption}[section]

\renewcommand{\footnote}{}
\numberwithin{equation}{section}

\newcommand{\R}{\mathbb{R}}

\newcommand{\diff}{{\,\rm{d}}}

\makeatletter
    
    \newcommand{\Rmnum}[1]{\expandafter\@slowromancap\romannumeral #1@}
\makeatother

\begin{document}
\begin{frontmatter}
      \title{Strong convergence rates of stochastic theta methods for index 1 stochastic differential algebraic equations under non-globally Lipschitz conditions\tnoteref{mytitlenote}}


      \tnotetext[mytitlenote]{This work was supported by National Natural Science Foundation of China (Nos. 12201552 and 11961029), Jiangxi Provincial Natural Science Foundation(No. 20242BAB23004), Yunnan Fundamental Research Projects (No. 202301AU070010) and Innovation Team of School of Mathematics and Statistics of Yunnan University (No. ST20210104).}

      \author{Lin Chen}
      \address{School of Statistics and Data Science, Jiangxi University of Finance and Economics, Nanchang, 330013, China}

      \author{Ziheng Chen\corref{mycorrespondingauthor}}
      \address{School of Mathematics and Statistics, Yunnan University, Kunming 650500, China}
      \cortext[mycorrespondingauthor]{Corresponding author. E-mail address: czh@ynu.edu.cn (Ziheng Chen).}

      \author{Jing Zhao}
      \address{School of Mathematics and Statistics, Yunnan University, Kunming 650500, China}

      \begin{abstract}
            This work investigates numerical approximations of index 1 stochastic differential algebraic equations (SDAEs) with non-constant singular matrices under non-global Lipschitz conditions. Analyzing the strong convergence rates of numerical solutions in this setting is highly nontrivial, due to both the singularity of the constraint matrix and the superlinear growth of the coefficients. To address these challenges, we develop an approach for establishing mean square convergence rates of numerical methods for SDAEs under global monotonicity conditions. Specifically, we prove that each stochastic theta method with $\theta \in [\frac{1}{2},1]$ achieves a mean square convergence rate of order $\frac{1}{2}$. Theoretical findings are further validated through a series of numerical experiments.
      \end{abstract}


\end{frontmatter}

\section{Introduction}


The dynamical behavior of random phenomenon, arising in fields such as physics, engineering, and finance, are typically described by stochastic differential equations (SDEs); see, e.g., \cite{friedman2006stochastic, mao2008stochastic, platen2010numerical}. However, when the states of a stochastic system are subject to additional constraints—such as conservation laws (e.g., Kirchhoff’s current laws in electrical networks) or positional restrictions (e.g., the motion of mass points confined to a surface)—the corresponding mathematical formulation must also incorporate algebraic equations to capture these constraints \cite{denk2008efficient, suthar2022explicitly}. Stochastic systems that simultaneously involve both algebraic constraints and differential equations are referred to as stochastic differential algebraic equations (SDAEs). Such equations enrich mathematical models by additionally incorporating static constraint mechanisms, thereby extending their applicability to a wide variety of real-world systems. At the same time, these constraints introduce strong coupling within the system, which significantly increases the complexity and challenge of their theoretical analysis. We refer to \cite{alabert2006linear, cong2014adjoint, nguyen2010stochastic, cong2012lyapunov, gliklikh2024stochastic, serea2025existence, the2023central, nguyen2024stochastic, thuan2024stability, winkler2003stochastic} and references therein for variety descriptions, properties and applications of SDAEs.

In this paper, we consider the following nonlinear SDAEs
\begin{equation}\label{eq:SDAE}
      A_t \text{d} X_t 
      = F(t,X_t) \text{d}t + G(t,X_t) \text{d} W_t,
      \quad t \in (0,T]
\end{equation}
with a random initial value $X_0$, where $F \colon [0,T] \times \mathbb{R}^{d} \to \mathbb{R}^{d}$ and 
$G \colon [0,T] \times \mathbb{R}^{d} \rightarrow \mathbb{R}^{d \times m}$ are Borel measurable functions. Here, $\{W_{t}\}_{t \in [0,T]}$ is an $m$-dimensional standard Brownian motion on the complete filtered probability space $(\Omega,\mathcal{F},\mathbb{P},\{\mathcal{F}_{t}\}_{t \in [0,T]})$ with the filtration $\{\mathcal{F}_{t}\}_{t \in [0,T]}$ satisfying the usual conditions. Besides, the matrix $A_{t} \in \mathbb{R}^{d \times d}$ is singular for each $t \in [0,T]$ and thus acts as the key feature to distinguish SDAEs \eqref{eq:SDAE} from the unconstrained SDEs. 
Moreover, we assume that \eqref{eq:SDAE} is an index 1 SDAEs system, i.e., the noise sources do not appear in the algebraic constraints and the constraints are globally uniquely solvable with respect to the algebraic variables \cite{serea2025existence}. Nowadays, SDAEs \eqref{eq:SDAE} have been widely used to model different evolution phenomena arising in various fields, such as constrained dynamics of protein folding problem \cite{neumaier1997molecular}, constrained Brownian motion in nanostructured materials \cite{raccis2011confined} and transient noise simulation of microelectronic circuits \cite{schein1998numerical}. Since the explicit solutions of most SDAEs are rarely available, numerical approximation has become a powerful tool to study the behaviour of their solutions. Up to now, much progress has been achieved in constructing and analyzing different kinds of numerical methods for SDAEs with constant $A_{t} \equiv A, t \in [0,T]$ under the global Lipschitz conditions and the linear growth conditions (see, e.g., \cite{avaji2019stability, kupper2009runge, kupper2012rungekutta, kupper2015stability, nair2021stochastic, schein1998numerical, sickenberger2009local, wang2011numerical, winkler2003stochastic}). However, such conditions are too restrictive in the sense that
some important SDAEs models, such as constrained stochastic Ginzburg--Landau equation and transient noise analysis model for circuits, only satisfy the local Lipschitz conditions with non-constant $A_{t}, t \in [0,T]$. As far as we know, \cite{tsafack2025pathwise} is the only work to cope with this situation recently, where the authors build a numerical scheme based on semi-implicit method for \eqref{eq:SDAE} and derive its pathwise convergence rate under non-global Lipschitz conditions. Given that strong convergence rates of numerical approximations play a crucial role in the design of efficient multilevel Monte Carlo methods \cite{giles2008improved}, the present work seeks to contribute to the study of mean square convergence rates of numerical approximations for \eqref{eq:SDAE} with global monotonicity conditions.

Taking an equidistant stepsize $\Delta = \frac{T}{K}$ for some $K \in \mathbb{N}$ and the uniform grids $t_{k} := k \Delta, k = 0,1,2, \cdots, K$, the stochastic theta method (STM) for each $\theta \in [\frac{1}{2},1]$ applied to \eqref{eq:SDAE} is given by $ x_0 = X_0$ and 
\begin{equation}\label{eq:ABEM}
      A_{t_{k+1}}x_{k+1}
      = 
      A_{t_k}x_k +\theta F(t_{k+1},x_{k+1}) \Delta 
      +
      (1-\theta) F(t_{k},x_{k}) \Delta 
      +
      G(t_k,x_k) \Delta W_{k}
\end{equation}
for $k = 0,1,2, \cdots, K-1$, where $x_{k}$ denotes the numerical approximation of $x_{t_{k}}$ and $\Delta{W_k} := W({t_{k+1}}) - W({t_k})$ is the Brownian increment. Benefiting the advantage of overcoming the difficulties caused by the superlinear growth of coefficients, STM has been successfully used in the strong convergence analysis for different kinds of unstrainted SDEs, e.g., \cite{chen2019stability, chen2025stochastic, jiang2020stationary, wang2020meansquare, zhang2022convergence}. However, the presence of algebraic constraints makes our situations more challenging, as the numerical method must not only approximate the differential part accurately but also preserve the constraint compatibility at each time step—a requirement that does not exist for unconstrained SDEs. 
To achieve our goal, we firstly transform \eqref{eq:SDAE} into an equivalent combination of an unconstrained SDEs \eqref{eq:inherent} and a nonlinear algebraic equation \eqref{eq:constraint} by means of the singular value decomposition of matrix and the implicit function theorem. With this preparation, we then apply the mathematical induction method to establish the well-posedness of numerical solutions $\{x_{k}\}_{k = 0,1,\cdots,K}$, including the existence of unique solution and the feasibility of algebraic constraints. Armed with the H\"{o}lder continuity of exact solution $\{X_t\}_{t \in [0,T]}$ and an upper error bound for the mean square error between $X_{t_{k}}$ and $x_{k}$, we finally prove that numerical solutions generated by each STM converge strongly with order $\frac{1}{2}$ in the mean square sense to exact solutions of \eqref{eq:SDAE} under the coupled monotonicity condition \eqref{asm:FG-FG} and the polynomial growth conditions \eqref{asm:FG-absF}--\eqref{asm:FG-absG}; see Theorem \ref{thm:convAt} for more details. Although Theorem \ref{thm:convAt} can directly yields the convergence rate for each STM of \eqref{eq:SDAE} with a constant matrix $A_{t} \equiv A,t \in [0,T]$, Theorem \ref{th:conv2A} establishes this result under more general and less restrictive assumptions.

Before concluding the introduction, it is worth noting that our results in Theorems \ref{thm:convAt} and \ref{th:conv2A} coincide with those of \cite{wang2020meansquare} in the case where the matrix $A_{t}$ is non-singular for each $t \in [0,T]$. In that work, the authors established the mean square convergence rates of STM for unconstrained SDEs under non-globally Lipschitz conditions. In addition to directly computing numerical solutions of SDAEs \eqref{eq:SDAE} as carried out in this work, one may also consider their equivalent equations \eqref{eq:constraint} and \eqref{eq:inherent} for numerical approximations, where \eqref{eq:inherent} can be directly analyzed using the results of \cite{andersson2017meansquare, hutzenthaler2012strong, tretyakov2013fundamental, wang2020meansquare}. However, solving the nonlinear algebraic equation \eqref{eq:constraint} requires iterative computations in practice and therteby substantially increases the overall computational cost, which does not need for our methods as the algebraic constraints have been embedded in the discrete method \eqref{eq:ABEM}. Finally, compared with the existing work \cite{tsafack2025pathwise}, this study develops a framework for establishing the strong convergence rates of SDAEs \eqref{eq:SDAE} under non-global Lipschitz conditions. This consitutes the main contribution in this work and provides a foundation that can be further extended to the error analysis of Milstein-type methods in our future research.

The remainder of this paper is organized as follows. The next section concerns assumptions and properties of the underlying SDAEs \eqref{eq:SDAE}. In Section \ref{sec:orderAt}, we establish the mean square convergence rate of each STM of \eqref{eq:SDAE} under non-global Lipschitz conditions. For the case of $A_{t} \equiv A,t \in [0,T]$ being a constant matrix, the convergence analysis is performed in Section \ref{sec:orderA}. Finally, several numerical experiments are given in Section \ref{sec:experiments} to illustrate these theoretical results.

\section{SDAEs of index 1}\label{SDAEs}
In this section, we formulate a set of assumptions for SDAEs \eqref{eq:SDAE} and establish their key properties through several lemmas, which will serve as the foundation for our subsequent convergence analysis. To this end, we start with some notation that will be used throughout the paper. For simplicity, the letter $C$ stands for a generic positive constant whose value may vary for each appearance, but independent of the stepsize $\Delta$ of the considered numerical method.  Let $\langle \cdot,\cdot \rangle$ and $|\cdot|$ be the Euclidean inner product and the corresponding norm of vectors in $\R^{d}$, respectively. By $B^{\top}$ we denote the transpose of a vector or a matrix $B$. For any matrix $B \in \R^{d \times m}$, we use $|B| := \sqrt{\rm{trace}(B^{\top}B)}$ to denote the trace norm of $B$, ${\rm{Im}}(B) := \{Bx: x \in \mathbb{R}^{m}\}$ the image of $B$, and ${\rm{Ker}}(B) := \{x \in \mathbb{R}^{m} : Bx = 0\}$ the kernel of $B$, respectively.

We first make the following assumption on the singular matrix $A_{t}$ for $t \in [0,T]$.

\begin{assumption}\label{asm:At}
      For any $t\in [0,T]$, the singular value decomposition of matrix $A_t$ takes the form
      $A_t = M \Sigma_t N$, where $M,N \in \mathbb{R}^{d \times d}$ are orthogonal matrices and $\Sigma_t= {\rm{diag}}(\sigma_{1}(t), \cdots, \sigma_{r}(t), 0, \cdots, 0)$ with $r \in \{1,2,\cdots,d\}$ is continuously differentiable. Moreover, there exist constants $\underline{\sigma}$ and $\overline{\sigma}$ such that for any $i \in \{1,2,\cdots,r\}$ and $t \in [0,T]$,
      \begin{align*}
            0 < \underline{\sigma}
            \leq
            \sigma_{i}(t)
            \leq
            \overline{\sigma} < \infty.
      \end{align*}
\end{assumption}

Given a matrix $B \in \mathbb{R}^{d \times d}$ (may be singular), the following definition states the definition of its pseudo inverse matrix (see, e.g., \cite{gentle2024matrix}), which is a generalization of the inverse matrix. Indeed, if $B$ is a non-singular matrix, then the pseudo inverse coincides with the inverse of the matrix $B$. 

\begin{definition}\label{def:pinverse}
      A matrix $B^- \in \mathbb{R}^{d \times d}$ is said to be the pseudo inverse of $B \in \mathbb{R}^{d \times d}$ if it satisfies all four conditions below:
      \begin{equation*}
            BB^{-}B = B, \quad  B^{-}BB^{-} = B^{-},
            \quad
            (BB^{-})^{\top} = BB^{-},
            \quad  
            (B^{-}B)^{\top} = B^{-}B.
      \end{equation*}
\end{definition}

Since the pseudo inverse always exists and is unique for any square matrix, one can denote $A_{t}^{-}$ the unique pseudo inverse matrix of $A_{t}$ for each $t \in [0,T]$. Then Assumption \ref{asm:At} implies $A_{t}^{-} = N^{\top} \Sigma_{t}^{-} M^{\top}$ with $\Sigma_{t}^{-} = {\rm{diag}}(\sigma_{1}^{-1}(t), \cdots, \sigma_{r}^{-1}(t), 0, \cdots, 0)$. Moreover, the time independent operator $$P := A_{t}^{-}A_{t} = N^{\top} \begin{pmatrix} I_r &  \\  & \mathbf{0}_{d-r} \end{pmatrix} N$$ is a projector along  ${\rm{Ker}}(A_{t})$, where $\mathbf{0}_{d-r}$ represents the $(d-r) \times (d-r)$ zero matrix. The projector $Q = I - P$ onto ${\rm{Ker}}(A_{t})$ is such that $A_{t}Q = 0$ and the projector $R = I - A_{t}A_{t}^{-}$ along ${\rm{Im}}(A_{t})$ satisfies $RA_{t} = 0$. By these special projectors, one can split the solution of SDAEs \eqref{eq:SDAE} into the differential and algebraic components as follows
\begin{equation}\label{eq:PQ}
      X_t = PX_t + QX_t := U_t+V_t,
      \quad
      U_t \in \text{Im}(P), \ V_t\in \text{Im}(Q).
\end{equation}
Since \eqref{eq:SDAE} is assumed to be of 1 index, 
we have $\text{Im}(G(t,x)) \subset \text{Im}(A_t)$ for $t \in [0,T], x \in \mathbb{R}^{d}$,
which implies
\begin{equation*}
      R G(t,x)=0, 
      \quad t \in [0,T],
      x \in \mathbb{R}^{d}. 
\end{equation*}
Together with $A_tV_t=A_tQX_t=0$ for $t \in [0,T]$, we obtain the nonlinear algebraic equation
\begin{equation}\label{eq:constraint}
      A_tV_t + R F(t,U_t+V_t)=0, \quad t \in [0,T].
\end{equation}

Denoting $J(t,x) := A_{t} + R F'_{x}(t,x), t \in [0,T], x \in \mathbb{R}^{d}$ the Jacobian matrix of \eqref{eq:constraint}, the assumption that \eqref{eq:SDAE} is of 1 index
indicates that $J(t,x)$ is nonsingular and the implicitly defined function exists globally and uniquely. Consequently, for any  $t\in [0,T]$ and $x\in\mathbb{R}^d$, there exists a unique solution (see, e.g., \cite{winkler2003stochastic}) denoted by $v=\hat{V}(t,x)$ of the following equation 
\begin{equation}\label{eq:v(t,x)}
      A_tv +R F(t,x+v)=0.
\end{equation}
It follows that for any $x \in \mathbb{R}^d$ and $t \in [0,T]$, 
\begin{equation}\label{def:Qx}
      x = Px + Qx, \quad \hat{V}(t,Px) = Qx + \hat{V}(t,x).
\end{equation}
Moreover, $V_t=\hat{V}(t,U_t), t \in [0,T]$ is the unique solution of \eqref{eq:constraint} and \eqref{eq:PQ} becomes
\begin{equation}\label{eq:XU}
      X_t = U_t + \hat{V}(t,U_t), \quad t\in [0,T].
\end{equation}
Multiplying \eqref{eq:SDAE} by $A_t^-$ leads to 
\begin{equation*}
      P \text{d} X_t 
      = A_t^- F(t,X_t) \text{d}t + A_t^- G(t,X_t) \text{d} W_t,
      \quad t \in [0,T].
\end{equation*}
which along with $U_t = PX_t, t \in [0,T]$ and \eqref{eq:XU} yields the following unconstrained SDEs
\begin{equation}\label{eq:inherent}
      U_{t} - U_{0}
      =
      \int_{0}^{t} A_{t}^{-}
      F(s,U_{s} + \hat{V}(s,U_{s})) \text{d}s
      +
      \int_{0}^{t} A_{t}^{-} G(s,U_{s} 
      + \hat{V}(s,U_s)) \text{d} W_s,
      \quad t \in [0,T]
\end{equation}
with the initial value $U_{0} = PX_{0}$ such that $RF(0,X_0)=0$. Consequently solving \eqref{eq:SDAE} is equivalent to solving the coupled system \eqref{eq:constraint} and \eqref{eq:inherent}. Besides, letting the measurable functions $f \colon [0,T] \times \mathbb{R}^{d} \to \mathbb{R}^d$ and $G \colon [0,T] \times \mathbb{R}^{d} \to \mathbb{R}^{d \times m}$ be given by 
$$f(t,x) := A_t^- F(t,x+\hat{V}(t,x)),
\quad
g(t,x) := A_t^- G(t,x+\hat{V}(t,x))$$
for any $t \in [0,T]$ and $x \in \mathbb{R}^{d}$, then \eqref{eq:inherent} can be rewritten as follows
\begin{equation}\label{eq:inherent_U}
      U_t-U_0
      =
      \int_0^t f(s,U_s) \text{d}s
      +
      \int_0^t g(s,U_s) \text{d} W_s,
      \quad t \in [0,T].
\end{equation}

To proceed, we put the following assumption on \eqref{eq:SDAE}.

\begin{assumption}\label{asm:FG}
      Suppose that the following conditions are satisfied.

      \noindent {\rm{(A1)}} There exist constants $L_1>0$, $\gamma\geq 1$ and $p_1 > 4\gamma-2$  such that for any $x,y \in\mathbb{R}^d$ and $t, s \in [0,T]$,
      \begin{align}\label{asm:FG-FG}
            2\big\langle Px-Py, 
            A_t^- F(t,x) - A_t^- F(t,y) \big\rangle 
            +
            (p_1-1)|A_t^- G(t,x) - A_t^- G(t,y))|^2
            \leq 
            L_1|x-y|^2,
      \end{align}
      and
      \begin{equation}\label{asm:FG-absF}
            |F(t,x)-F(s,y)|
            \leq
            C\big((1+|x|+|y|)^{\gamma-1} |x-y|
            + (1+|x|+|y|)^{\gamma} |t-s|\big),
      \end{equation}
      as well as
      \begin{equation}\label{asm:FG-absG}
            |G(t,x)-G(s,y)|^{2}
            \leq
            C\big((1+|x|+|y|)^{\gamma-1} |x-y|^{2}
            + (1+|x|+|y|)^{\gamma+1} |t-s|^2\big).
      \end{equation}

      \noindent {\rm{(A2)}} For all $t \in [0,T], x \in \mathbb{R}^{d}$, the Jacobian matrix $J(t,x):= A_t + RF'_x(t,x)$ is continuous and possesses a globally bounded inverse, i.e., there exists a positive constant $L_J$  such that
      \begin{equation*}
            |J(t,x)^{-1}| \leq L_{J},\quad t \in [0,T], x \in \mathbb{R}^{d}.
      \end{equation*}

      \noindent {\rm{(A3)}} The initial value $X_0$ satisfies $RF(0,X_0) = 0$ almost surely.
\end{assumption}


Under Assumptions \ref{asm:At} and \ref{asm:FG}, \eqref{eq:inherent} admits a unique solution $\{U_{t}\}_{t \in [0,T]}$ \cite{serea2025existence}, which in combination with \eqref{eq:XU} guarantees the existence of unique solution $\{X_{t}\}_{t \in [0,T]}$ to \eqref{eq:SDAE}. Before ending this section, we collect several properties related to the underlying equations. The first lemma concerns the uniform estimate on $A_{t}, t \in [0,T]$ as well as its pseudo inverse matrix and derivetive matrix.

\begin{lemma}\label{lem:prop:A}
      Let $A_t' := \frac{\text{d}}{\text{d}t} A_t =M \Sigma_t' N$ for $t \in [0,T]$ and suppose that Assumption \ref{asm:At} holds. Then there exists $C > 0$ such that for any $t \in [0,T]$,
      \begin{equation*}
            |A_t|
            \leq
            \sqrt{r} d \overline{\sigma},
            \quad 
            |A_t^-|
            \leq
            \frac{\sqrt{r}d}{\underline{\sigma}},
            \quad
            |A_{t}'|
            \leq 
            C.
      \end{equation*}
\end{lemma}

\textbf{Proof.}\quad
      Owing to Assumption \ref{asm:At}, we have that for any $t \in [0,T]$,
      \begin{equation*}
            |\Sigma_t|^{2} 
            = \sum_{i=1}^{r} (\sigma_i(t))^{2}
            \leq 
            r \overline{\sigma}^{2},
            \quad
            |\Sigma_t^-|^{2}
            =
            \sum_{i=1}^{r} (\sigma_i^{-1}(t))^{2}
            \leq 
            \frac{r}{\underline{\sigma}^{2}}.
      \end{equation*}
      Since $M,N$ are $d$-dimensional orthogonal matrix, we use $|M| = |N| = \sqrt{d}$ to get
      \begin{equation*}
             |A_t|
             \leq
             |M|| \Sigma_t| |N|
             \leq \sqrt{r} d \overline{\sigma},
             \quad 
             |A_t^-|
             \leq
             |M|| \Sigma_t^-| |N|
             \leq 
             \frac{\sqrt{r}d}{\underline{\sigma}},
      \end{equation*}
      and
      \begin{equation*}
            |A_t'|
            \leq
            |M|| \Sigma_t'| |N|
            \leq 
            d\sqrt{\sum_{i=1}^r (\sigma_i'(t))^{2}}
            \leq 
            d\sqrt{r}
            \sqrt{\max\Big\{\Big(\max_{t \in [0,T]}
            |\sigma_i'(t)|\Big)^{2}: 
            i = 1,2,\cdots,r\Big\}}
      \end{equation*}
      due to the continuity of $\Sigma_i'(t)$ on $[0,T]$ for $i = 1,2,\cdots,r$. Thus we complete the proof.\hfill$\square$

The next lemma provides different properties of function $V(t,x)$ determined by \eqref{eq:v(t,x)}.

\begin{lemma}\label{lem-v}
      Denote $\hat{L} := \Big(\sup\limits_{t \in [0,T]}|A_t|\Big)L_{J} + d$ and let the implicit function $\hat{V} \colon [0,T] \times \mathbb{R}^d \rightarrow \mathbb{R}^d, (t,x) \mapsto V(t,x)$ be the unique solution of \eqref{eq:v(t,x)}.
      Suppose that Assumption \ref{asm:At} and \ref{asm:FG} hold. Then there exists a constant $C > 0$ such that for any $s,t\in[0,T]$ and $u,v \in \mathbb{R}^d$, 
      \begin{gather}
            \label{lem-v-rs1}
            \big|\hat{V}(t,u)\big|
            \leq
            C\big(1 + |u|\big),
            \quad
            \big|u + \hat{V}(t,u)\big|
            \leq 
            C\big(1 + |Pu|\big),
            \\\label{lem-v-rs2}
            \big|\hat{V}(t,u) - \hat{V}(s,v)\big|
            \leq
            \hat{L}\big|u-v\big| 
            + C\big(1 + |u| + |v|\big)^{\gamma}
            \big|t-s\big|,
            \\\label{lem-v-rs3}
            \big|\big(u + \hat{V}(t,u)\big)
            - \big(v + \hat{V}(s,v)\big)\big|
            \leq
            \big(1 + \hat{L}\big)\big|Pu - Pv\big| 
            + C\big|t-s\big|.
      \end{gather}
\end{lemma}

\textbf{Proof.}\quad
      For function $h(t,u,v) := A_t v +RF(t,u+v), t \in [0,T], u,v \in \mathbb{R}^{d}$, the definition of $\hat{V}$ in \eqref{eq:v(t,x)} gives that 
      \begin{equation}\label{eq:373737}
            h\big(t,u,\hat{V}(t,u)\big) \equiv 0,
            \quad t \in [0,T], u \in \mathbb{R}^{d}.
      \end{equation}
      Then for any $t \in [0,T]$ and $u \in \mathbb{R}^{d}$, 
      \begin{align*}
            0 
            =
            \frac{\partial}{\partial u} 
            h\big(t,u,\hat{V}(t,u)\big)
            = 
            h'_u \big(t,u,\hat{V}(t,u) \big) 
            +
            \hat{V}'_u(t,u) h'_v 
            \big(t,u,\hat{V}(t,u)\big),
      \end{align*}
      which together with {\rm{(A2)}} in Assumption \ref{asm:FG} implies 
      \begin{align*}
            \big|\hat{V}'_{u}(t,u)\big| 
            =&~
            \left| -h'_{u} \big(t,u,\hat{V}(t,u)\big) 
            \big(h'_{v}\big(t,u,\hat{V}(t,u)\big)
            \big)^{-1}\right| \nonumber 
            \\=&~
            \left|-R F'_x(t,u+\hat{V}(t,u)) 
            \big( A_t +R F'_x(t,u+\hat{V}(t,u))
            \big)^{-1}\right| \nonumber 
            \\=&~
            \left|\big(A_{t} - \big(A_{t}
            + R F'_x(t,u+\hat{V}(t,u))\big) \big)
            \big(A_t +R F'_x(t,u+\hat{V}(t,u))\big)^{-1}
            \right| \nonumber 
            \\=&~
            \left| A_t J(t,u+\hat{V}(t,u))^{-1} 
            - I_d \right| \nonumber 
            \\\leq&~
            \hat{L}.
      \end{align*}
      As a consequence of the mean value theorem, one gets that for any $u,v \in \mathbb{R}^d$ and $t \in [0,T]$,
      \begin{align}\label{lem-v-uv}
            \big|\hat{V}(t,u) - \hat{V}(t,v)\big|
            =
            \big|\hat{V}_{u}'(t,v+\theta(u-v))(u-v)\big|
            \leq
            \hat{L}\big|u-v\big|,
      \end{align}
      which results in
      \begin{align}\label{lem-v-absv}
            \big|\hat{V}(t,u)\big|
            \leq
            \hat{L}\big|u\big| + \big|\hat{V}(t,0)\big| 
            \leq
            \hat{L}\big|u\big| 
            + \sup_{t \in [0,T]}\big|\hat{V}(t,0)\big|,
            \quad u \in \R^{d}, t \in [0,T],
      \end{align}
      i.e., the first inequality of \eqref{lem-v-rs1} holds. Owing to \eqref{def:Qx}, we have that for any $u \in \mathbb{R}^{d}$ and $t \in [0,T]$,
      \begin{align*}
            \big|u + \hat{V}(t,u)\big|
            =
            \big|Pu + \hat{V}(t,Pu)\big|
            \leq
            \big|Pu\big| + \big|\hat{V}(t,Pu)\big|
            \leq
            C\big(1 + |Pu|\big),
      \end{align*}
      which gives the second inequality of \eqref{lem-v-rs1}.

      On the one hand, \eqref{eq:373737} shows that for any $t \in [0,T]$ and $u \in \mathbb{R}^{d}$, 
      \begin{align*}
            0 
            =&~
            \frac{\partial}{\partial t} 
            h\big(t,u,\hat{V}(t,u)\big) \nonumber 
            \\=&~
            A'_t \hat{V}(t,u) + A_t \hat{V}'_t(t,u) 
            +
            R\big( F'_t \big(t,u+\hat{V}(t,u)\big) 
            +
            F'_x\big(t,u+\hat{V}(t,u)\big) 
            \hat{V}'_t(t,u) \big),
      \end{align*}
      which along with {\rm{(A2)}} in Assumption \ref{asm:FG}, \eqref{lem-v-absv} and $|F_{t}'(t,x)| \leq C\big(1 + |x|\big)^{\gamma}$ for any $t \in [0,T], x \in \R^{d}$ leads to
      \begin{align*}
            \big|\hat{V}'_t(t,u)\big| 
            =&~
            \left|\big(A'_t \hat{V}(t,u) 
            + R F'_t \big(t,u+\hat{V}(t,u)\big) \big)
            \big( A_t + RF'_x(t,u+\hat{V}(t,u)) \big)^{-1}\right| \nonumber 
            \\\leq&~
            L_{J}\big(|A'_t| |\hat{V}(t,u)|
            + |R| |F'_t(t,u+\hat{V}(t,u))|\big) \nonumber
            \\\leq&~
            L_{J}\big(|A'_t| |\hat{V}(t,u)|
            + C|R|(1 + |u+\hat{V}(t,u)|)^{\gamma}\big) \nonumber
            \\\leq&~
            C\big(1+|u|^{\gamma}\big).
      \end{align*}
      Combining with \eqref{lem-v-uv} yields that for any $u,v \in \mathbb{R}^{d}$ and $t,s \in [0,T]$,
      \begin{align*}
            \big|\hat{V}(t,u)-\hat{V}(s,v)\big|
            \leq&~ 
            \big|\hat{V}(t,u)-\hat{V}(t,v)\big| 
            +
            \big|\hat{V}(t,v)-\hat{V}(s,v)\big| \nonumber 
            \\\leq&~ 
            \hat{L}\big|u-v\big| 
            + C\big(1+|v|\big)^{\gamma}
            \big|t-s\big| \nonumber
            \\\leq&~ 
            \hat{L}\big|u-v\big| 
            + C\big(1 + |u| + |v|\big)^{\gamma}
            \big|t-s\big|,
      \end{align*}
      i.e., \eqref{lem-v-rs2} holds. On the other hand, \eqref{def:Qx} and \eqref{lem-v-rs1} indicate that for any $u,v \in \mathbb{R}^{d}$ and $t,s \in [0,T]$,
      \begin{align*}
            \big|\big(u+\hat{V}(t,u)\big)
            - \big(v+\hat{V}(s,v)\big)\big|
            =&~ \nonumber
            \big|\big(Pu + \hat{V}(t,Pu)\big)
            - \big(Pv + \hat{V}(s,Pv)\big)\big|
            \\\leq&~
            \big|Pu - Pv\big|  \nonumber
            +
            \big|\hat{V}(t,Pu) -\hat{V}(s,Pv)\big|
            \\\leq&~
            \big(1 + \hat{L}\big)\big|Pu - Pv\big| 
            + C\big(1 + |Pu| + |Pv|\big)^{\gamma}|t-s|.
      \end{align*}
      which means that \eqref{lem-v-rs3} holds. Thus we complete the proof.\hfill$\square$

To develop the boundness of higher moment and the H\"{o}lder continuity of solution process $\{U_{t}\}_{t \in [0,T]}$, Lemmas \ref{lem-f} and \ref{lem-v-3} below focus on the coefficients of \eqref{eq:inherent_U}.

\begin{lemma}\label{lem-f}
      Suppose that Assumptions \ref{asm:At} and \ref{asm:FG} hold. Then for any $x,y \in \mathbb{R}^{d}$ and $t \in [0,T]$, 
      \begin{equation*}
            f(t,x) = f(t,Px),
            \quad 
            g(t,x)=g(t,Px),
            \quad
            \big\langle x,f(t,y) \big\rangle 
            = 
            \big\langle Px,f(t,y) \big\rangle.
      \end{equation*}
\end{lemma}

\textbf{Proof.}\quad
      For any $x,y \in \mathbb{R}^{d}$ and $t \in [0,T]$, the expressions $f(t,x) := A_t^- F\big(t,x+\hat{V}(t,x)\big), g(t,x) := A_t^- G\big(t,x+\hat{V}(t,x)\big)$ and \eqref{def:Qx} enable us to get
      \begin{align*}
            f(t,x)
            =&~
            A_t^- F\big(t,Px+Qx+\hat{V}(t,x)\big)
            \nonumber 
            =
            A_t^- F\big(t,Px+\hat{V}(t,Px)\big)
            =
            f(t,Px),
      \end{align*}
      and
      \begin{align*}
            g(t,x)
            =&~
            A_t^- G\big(t,Px+Qx+\hat{V}(t,x)\big)
            \nonumber 
            =
            A_t^- G\big(t,Px+\hat{V}(t,Px)\big)
            =
            g(t,Px).
      \end{align*}
      Noting that $P = A_t^-A_t$ is symmetric and satisfies $PA_t^-=A_t^-$, we further derive
      \begin{align*}
            \big\langle x,f(t,y) \big\rangle
            =&~ \nonumber
            \big\langle x,A_t^- 
            F\big(t,y+\hat{V}(t,y)\big) \big\rangle
            \\=&~
            \big\langle x, PA_t^- 
            F(t,y+\hat{V}(t,y)) \big\rangle \nonumber 
            \\=&~
            \big\langle Px, A_t^- 
            F(t,y+\hat{V}(t,y)) \big\rangle \nonumber 
            \\=&~
            \big\langle Px,f(t,y) \big\rangle,
      \end{align*}
      which gives the desired result.\hfill$\square$

\begin{lemma}\label{lem-v-3}
      Suppose that Assumptions \ref{asm:At} and \ref{asm:FG} hold. Then  for any $x,y \in \mathbb{R}^{d}$ and $0 \leq t \leq T$, \begin{align}\label{lem:fg:rs1}
            2\big\langle x - y, 
            f(t,x) - f(t,y) \big\rangle 
            + \big(p_{1}-1\big)
            \big|g(t,x) - g(t,y)\big|^{2}
            \leq
            2L_{1}\big(1 + \hat{L}^{2}\big)
            \big|x - y\big|^{2}.
      \end{align}
      Besides, for any $2 \leq p < p_{1}$, there exists a constant $C > 0$ such that for any $x \in \R^{d}$ and $t \in [0,T]$, 
      \begin{equation}\label{lem:fg:rs4}
            2\big\langle x, f(t,x) \big\rangle 
            + \big(p-1\big)\big|g(t,x)\big|^{2}
            \leq 
            C\big(1 + |x|^{2}\big).
      \end{equation}
\end{lemma}

\textbf{Proof.}\quad
       For any $x,y \in \mathbb{R}^{d}$ and $0 \leq t \leq T$, let $\hat{V}(t,x)$ and $\hat{V}(t,y)$ be the unique solutions of equations $A_tv + RF(t,x+v) = 0$ and $A_tv+RF(t,y+v)=0$, respectively.  
       Then $P\hat{V}(t,x) = A_t^-A_t\hat{V}(t,x) = 0$ and Lemma \ref{lem-f} ensure that
       \begin{align*}
             &~2\big\langle x-y,f(t,x) -f(t,y) \big\rangle 
             + 
             \big(p_{1} - 1\big)
             \big|g(t,x)-g(t,y)\big|^{2}
             \\=&~
             2\big\langle P(x-y),
             f(t,x) -f(t,y) \big\rangle 
             + \big(p_{1} - 1\big)
             \big|g(t,x)-g(t,y)\big|^{2}
             \\=&~
             2\big\langle P\big(x+\hat{V}(t,x)\big)
             - P\big(y+\hat{V}(t,y)\big), 
             A_t^- F\big(t,x+\hat{V}(t,x)\big) 
             - A_t^- F\big(t,y+\hat{V}(t,y)\big)
             \big\rangle
             \\&~+
             \big(p_{1} - 1\big)
             \big|A_t^-G\big(t,x+\hat{V}(t,x)\big) 
             - A_t^-G\big(t,y+\hat{V}(t,y)\big)\big|^{2}
             \\\leq&~
             L_{1}\big|\big(x + \hat{V}(t,x)\big)
             - \big(y + \hat{V}(t,y)\big)\big|^{2}
             \\\leq&~
             2L_{1}\big|x-y\big|^{2} 
             + 
             2L_{1}\big|\hat{V}(t,x)-\hat{V}(t,y)\big|^{2}
             \\\leq&~
             2L_{1}\big(1 + \hat{L}^2\big)
             \big|x-y\big|^{2},
      \end{align*}
      where Assumption \ref{asm:FG} and \eqref{lem-v-uv} have been used. 
      Noting that $2 \leq p < p_{1}$, the Young inequality implies 
      \begin{align*}
            \big(p-1\big)\big|g(t,x)\big|^{2}
            =&~
            \big(p-1\big)\big|g(t,x) 
            - g(t,0) + g(t,0)\big|^{2}
            \nonumber 
            \\\leq&~
            \big(p-1\big)\bigg(
            \Big(1+\frac{p_1-p}{p-1}\Big)
            \big|g(t,x) - g(t,0)\big|^{2} 
            +
            \Big(1+\frac{p-1}{p_1-p}\Big)
            \big|g(t,0)\big|^{2}\bigg)
            \nonumber 
            \\\leq&~
            \big(p_{1}-1\big)
            \big|g(t,x) - g(t,0)\big|^{2} 
            + 
            \frac{(p-1)(p_{1}-1)}{p_{1}-p}\bigg(
            \sup_{t \in [0,T]}
            \big|g(t,0)\big|^{2}\bigg).
      \end{align*}
      It follows from \eqref{lem:fg:rs1} with $y = 0$ that 
      \begin{align*}
            &~2\big\langle x,f(t,x) \big\rangle 
            + \big(p-1\big)\big|g(t,x)\big|^{2}
            \nonumber 
            \\\leq&~
            2\big\langle x,f(t,x) -f(t,0) \big\rangle 
            + 
            \big(p_{1} - 1\big)
            \big|g(t,x)-g(t,0)\big|^{2}
            + 2\big\langle x,f(t,0) \big\rangle + C
            \nonumber
            \\\leq&~
            2L_{1}\big(1 + \hat{L}^{2}\big)\big|x\big|^{2}
            + \big|x\big|^{2} + \big|f(t,0)\big|^{2} + C
            \nonumber
            \\\leq&~
            C\big(1 + \big|x\big|^{2}\big).
      \end{align*}
      Then we finish the proof.\hfill$\square$

Based on the above preparations, we are able to present the moment estimates for solution processes $\{U_{t}\}_{t \in [0,T]}$ and $\{X_{t}\}_{t \in [0,T]}$.

\begin{lemma}\label{lem:bound}
      Suppose that Assumptions \ref{asm:At} and \ref{asm:FG} hold. Then for any $2 \leq p < p_{1}$, there exists a constant $C > 0$ such that 
      \begin{align*}
            \sup_{0 \leq t \leq T} 
            \mathbb{E}\big[|U_{t}|^{p}\big] \leq C,
            \quad
            \sup_{0 \leq t \leq T} 
            \mathbb{E}\big[|X_{t}|^{p}\big] \leq C.
      \end{align*}
\end{lemma}

\textbf{Proof.}\quad
      Applying the It\^{o} formula to \eqref{eq:inherent_U} indicates that for any $t \in [0,T]$,
      \begin{align*}
            \mathbb{E}\big[|U_{t}|^{p}\big]
            =&~
            \mathbb{E}\big[|U_{0}|^{p}\big]
            + 
            \mathbb{E}\bigg[\int_{0}^{t} 
            \Big(p|U_s|^{p-2} \big\langle U_{s}, 
            f(s,U_s) \big\rangle 
            +
            \frac{p}{2} \big|U_{s}\big|^{p-2} 
            \big|g(s,U_{s})\big|^{2}
            \\&~
            +\frac{p(p-2)}{2} \big|U_{s}\big|^{p-4} 
            \big|(U_s)^{\top} g(s,U_s)\big|^{2}
            \Big)\text{d}s \bigg] 
            \\\leq&~
            \mathbb{E}\big[|U_{0}|^{p}\big]
            + 
            \frac{p}{2} \mathbb{E}\bigg[
            \int_{0}^{t} \big|U_{s}\big|^{p-2} 
            \big(2\big\langle U_{s}, 
            f(s,U_{s}) \big\rangle 
            + (p-1)\big|g(s,U_{s})\big|^{2}
            \big)\text{d}s \bigg]
            \\\leq&~
            \mathbb{E}\big[|U_{0}|^{p}\big]
            + 
            \frac{p}{2}C \mathbb{E}\bigg[
            \int_{0}^{t} \big|U_{s}\big|^{p-2}
            \big(1 + \big|U_{s}\big|^{2}\big)
            \text{d}s \bigg]
            \\\leq&~
            \big(CT + \mathbb{E}\big[|U_{0}|^{p}\big]\big)
            + 
            C \int_{0}^{t} \mathbb{E}
            \big[|U_{s}|^{p}\big]\text{d}s
      \end{align*}
      due to \eqref{lem:fg:rs4} and the Young inequality. By the Gronwall inequality, we obtain  
      \begin{align}\label{eq:lem:fg:}
            \mathbb{E}\big[|U_{t}|^{p}\big]
            \leq
            \big(CT + \mathbb{E}\big[|U_{0}|^{p}\big]\big)
            \exp\bigg(C \int_{0}^{t} \text{d}s\bigg) 
            \leq
            \big(CT + \mathbb{E}
            \big[|U_{0}|^{p}\big]\big)e^{CT},
            \quad t \in [0,T].
      \end{align}
      From \eqref{eq:XU} and \eqref{lem-v-rs1}, it follows that for any $t \in [0,T]$,
      \begin{align*}
            \mathbb{E}\big[|X_{t}|^{p}\big]
            =
            \mathbb{E}\big[|U_t 
            + \hat{V}(t,U_t)|^{p}\big]
            \leq
            2^{p}\big(\mathbb{E}\big[|U_t|^{p}\big]
            + 
            \mathbb{E}\big[|\hat{V}(t,U_t)|^{p}\big]\big)
            \leq
            C\big(1 + \mathbb{E}\big[|U_t|^{p}\big]\big),
      \end{align*}
      which together with \eqref{eq:lem:fg:} completes the proof.\hfill$\square$

The last lemma provides the H\"{o}lder continuity of $\{U_{t}\}_{t \in [0,T]}$ and $\{X_{t}\}_{t \in [0,T]}$. 

\begin{lemma}\label{lem:UX}
      Suppose that Assumption \ref{asm:At} and \ref{asm:FG} hold. Then for any $2\leq q<p_1/\gamma$, there exists a constat $C > 0$ such that for any $s,t \in [0,T]$, 
      \begin{align*}
            \mathbb{E}\big[|U_t - U_s|^{q}\big]
            +
            \mathbb{E}\big[\big|X_t
            - X_s\big|^{q}\big]
            \leq 
            C\Big(\big(t-s\big)^{q}
            +
            \big(t-s\big)^{\frac{q}{2}}\Big).
      \end{align*}
\end{lemma}

\textbf{Proof.}\quad
      Without loss of generality we assume $s < t$. 
      Since $RF(\rho,X_{\rho}) = 0$ for $\rho \in [s,t]$, we have $X_{\rho} \in \mathcal{M}_{\rho}$, which along with $PX_{\rho} = U_{\rho}$ and \eqref{def:Qx=v} implies $X_{\rho} = PX_{\rho} + \hat{V}(\rho, PX_{\rho}) =  U_{\rho} + \hat{V}(\rho,U_{\rho})$. Applying \eqref{asm:FG-absF}, \eqref{lem-v-rs1} and Lemma \ref{lem:prop:A} leads to
      \begin{align*}
            \big|f(\rho,U_{\rho})\big|
            =&~
            \big|A_{\rho}^-F(\rho,U_{\rho} 
            + \hat{V}(\rho,U_{\rho}))\big|
            =
            \big|A_{\rho}^-F(\rho,X_{\rho})\big|
            \\\leq&~
            \big|A_{\rho}^-\big|
            \big(\big|F(\rho,X_{\rho}) - F(\rho,0)\big|
            + \big|F(\rho,0)\big|\big)
            \\\leq&~
            \Big(\sup_{\rho \in [0,T]}
            \big|A_{\rho}^-\big|\Big)
            \bigg(C\big(1 +
            \big|X_{\rho}\big|\big)^{\gamma}
            + 
            \Big(\sup_{\rho \in [0,T]}
            \big|F(\rho,0)\big|\Big)\bigg)
            \\\leq&~
            C\big(1 +
            \big|X_{\rho}\big|\big)^{\gamma}
            =
            C\big(1 + \big|U_{\rho} + \hat{V}
            (\rho,U_{\rho})\big|\big)^{\gamma}
            \\\leq&~
            C\big(1 + \big|U_{\rho}\big|^{\gamma}\big).
      \end{align*}
      Similarly, we utilize to get \eqref{asm:FG-absG}, \eqref{lem-v-rs1} and Lemma \ref{lem:prop:A} to derive
      \begin{align*}
            \big|g(\rho,U_{\rho})\big|^{2}
            =&~
            \big|A_{\rho}^- G(\rho,X_{\rho})\big|^{2}
            \leq
            2\big|A_{\rho}^-\big|^{2}\Big( \big|G(\rho,X_{\rho}) - G(\rho,0)\big|^{2}
            + \big|G(\rho,0)\big|^{2}\Big)
            \\\leq&~
            2\Big(\sup_{\rho \in [0,T]}
            \big|A_{\rho}^-\big|\Big)^{2}
            \bigg(C\big(1 + \big|X_{\rho}\big|\big)^{\gamma+1}
            + \Big(\sup_{\rho \in [0,T]}
            \big|G(\rho,0)\big|\Big)^{2}\bigg)
            \\\leq&~
            C\big(1 + \big|U_{\rho}\big|^{\gamma+1}\big).
      \end{align*}
      It follows from \eqref{eq:inherent_U}, the H\"older inequlity and the moment inequality for stochastic integrals (see, e.g., \cite[Theorem 7.1 in Chapter 1]{mao2008stochastic}) that
      \begin{align*}
            &~\mathbb{E}\big[|U_t - U_s|^{q}\big]
            \leq
            2^{q}\mathbb{E}\bigg[\Big|
            \int_s^t f({\rho},U_{\rho}) \text{d}{\rho}
            \Big|^{q}\bigg]
            +
            2^{q}\mathbb{E}\bigg[\Big|
            \int_s^t g({\rho},U_{\rho}) \text{d}W_{\rho}
            \Big|^{q}\bigg]
            \\\leq&~
            C\big(t-s\big)^{q-1}\int_s^t             \mathbb{E}\big[\big|f({\rho},U_{\rho}) 
            \big|^{q}\big]\text{d}{\rho}
            +
            C\big(t-s\big)^{\frac{q}{2}-1}
            \int_s^t \mathbb{E}\big[\big|
            g({\rho},U_{\rho})\big|^{q}\big] 
            \text{d}{\rho}
            \\\leq&~
            C\big(t-s\big)^{q-1}\int_s^t 
            \Big( 1 + \mathbb{E}\big[
            \big|U_{\rho}\big|^{q\gamma}\big]
            \Big)\text{d}{\rho}
            +
            C\big(t-s\big)^{\frac{q}{2}-1}
            \int_s^t \Big( 1 + \mathbb{E}\big[
            \big|U_{\rho}\big|^{\frac{q(\gamma+1)}{2}}\big]
            \Big)\text{d}{\rho}
            \\\leq&~
            C\Big(\big(t-s\big)^{q}
            +
            \big(t-s\big)^{\frac{q}{2}}\Big)       
      \end{align*}
      due to Lemma \ref{lem:bound} and the fact $\frac{q(\gamma+1)}{2} \leq q\gamma < p_{1}$. In combination with \eqref{eq:XU}, \eqref{lem-v-rs2} and Lemma \ref{lem:bound}, we have
      \begin{align*}
            \mathbb{E}\big[\big|X_t
            - X_s\big|^{q}\big]
            = &~
            \mathbb{E}\big[\big|\big(U_t 
            + \hat{V}(t,U_t)\big) - \big(U_s 
            + \hat{V}(s,U_s)\big)\big|^{q}\big]
            \\\leq&~
            2^{q}\mathbb{E}\big[\big|U_t 
            - U_s \big|^{q}\big]
            +
            2^{q}\mathbb{E}\big[\big|\hat{V}(t,U_t)
            - \hat{V}(s,U_s)\big|^{q}\big]
            \\\leq&~
            2^{q}\mathbb{E}\big[\big|U_t 
            - U_s \big|^{q}\big]
            +
            2^{q}\Big(\hat{L}\mathbb{E}\big[\big|U_t 
            - U_s \big|^{q}\big]
            \\&~+
            C\mathbb{E}\big[\big(1 + \big|U_t\big|
            + \big|U_s\big|\big)^{q\gamma}\big]
            \big(t-s\big)^{q}\Big)
            \\\leq&~
            C\Big(\big(t-s\big)^{q}
            +
            \big(t-s\big)^{\frac{q}{2}}\Big).
      \end{align*}
      Thus we complete the proof.\hfill$\square$

\section{Mean square convergence rates of STM with non-constant $A_t $}\label{sec:orderAt}
This section aims to establish the mean square convergence rates of each STM for SDAEs \eqref{eq:SDAE} under non-global Lipschitz conditions. To achieve this, we first face with the question of well-posedness of numerical solutions $\{x_{k}\}_{k = 0,1,\cdots,K}$
because of the implicitness of STMs \eqref{eq:ABEM}. Different from the implicit methods for unconstrained SDEs, the feasibility of algebraic constraints for numerical solutions of SDAEs must be taken into consideration. For this purpose, we define the constraint sets
\begin{equation}\label{def:mt}
      \mathcal{M}_t
      =
      \{ x\in\mathbb{R}^{d} : RF(t,x)=0 \},
      \quad t \in [0,T].
\end{equation}
and obtain that for any $t \in [0,T]$ and $x \in \mathcal{M}_t$, 
\begin{gather}\label{def:Qx=v}
      x = Px+\hat{V}(t,Px),
      \quad Qx = \hat{V}(t,Px),
      \\\label{def:Ff}
      A_t^- F(t,x) = f(t,Px),
      \quad
      A_t^- G(t,x) = g(t,Px).
\end{gather}
Let us prove the existence of unique numerical solutions $\{x_{k}\}_{k = 0,1,\cdots,K}$ satisfying $x_k \in \mathcal{M}_{t_k}$ for any $k = 0,1,2, \cdots, K$.

\begin{lemma}\label{lem:def}
      Suppose that Assumptions \ref{asm:At} and \ref{asm:FG} hold and let $\theta \in [\frac12,1]$. If $\Delta < \frac{1}{\theta{L_1}(1 + \hat{L}^{2})}$, then the ATM method \eqref{eq:ABEM} is well defined and $x_k \in \mathcal{M}_{t_k}$ for any $k = 0,1,2, \cdots, K$.
\end{lemma}

\textbf{Proof.}\quad
      The proof is by mathematical induction on $k$. Owing to {\rm{(A3)}} in Assumption \ref{asm:FG}, we have $x_0 \in \mathcal{M}_{t_0}$. Assume $x_k \in \mathcal{M}_{t_k}$ and define 
      \begin{align*}
            H(x)
            := 
            A_{t_{k+1}}x - \theta F(t_{k+1},x)\Delta 
            - \Big(A_{t_k}x_k + 
            (1-\theta) F(t_{k},x_{k})\Delta 
            + G(t_k,x_k) \Delta W_k\Big),
            \quad x \in \R^{d},
      \end{align*}
      which enables us to rewrite \eqref{eq:ABEM} as $H(x_{k+1})=0$. Because of $R = I-A_{t_{k+1}} A_{t_{k+1}}^-$, $H(x)=0$ is equivalent to
      \begin{align*}
      \begin{cases}
            A_{t_{k+1}}A_{t_{k+1}}^- H(x) = 0,
            \\
            RH(x)= 0.
      \end{cases}
      \end{align*}
      Noting the equivalence between $A_{t_{k+1}}A_{t_{k+1}}^- H(x)= 0$ and  $A_{t_{k+1}}^- H(x)= 0$ due to $A_{t_{k+1}}^-A_{t_{k+1}}A_{t_{k+1}}^-=A_{t_{k+1}}^-$, we further deduce that $H(x) = 0$ is equivalent to 
      \begin{align}\label{lem:def:A^-H}
      \begin{cases}
            A_{t_{k+1}}^- H(x)= 0,
            \\
            R H(x)= 0.
      \end{cases}
      \end{align}
      Now we claim that \eqref{lem:def:A^-H} admits a unique solution $x = x_{k+1}$. In fact, employing \eqref{lem:fg:rs1}, 
      $\Delta < \frac{1}{\theta{L_1}(1 + \hat{L}^{2})}$ and the uniform monotonicity theorem in \cite[Theorem C.2]{stuart1996dynamical} indicates that the following equation
      \begin{align*}
            x - \theta f(t_{k+1},x) \Delta 
            - \big(A_{t_{k+1}}A_{t_k}^-x_k
            + (1-\theta)A_{t_{k+1}}^-F(t_{k},x_{k}) \Delta
            + A_{t_{k+1}}^-G(t_{k},x_{k}) \Delta W_k\big)
            = 0, \quad x \in \mathbb{R}^{d}
      \end{align*}
      possesses a unique solution, denoted by $u_{k+1} \in \mathbb{R}^{d}$, i.e., 
      \begin{align}\label{lem:eq:u}
            u_{k+1} - \theta f(t_{k+1},u_{k+1}) \Delta 
            - \big(A_{t_{k+1}}A_{t_k}^-x_k
            + (1-\theta)A_{t_{k+1}}^-F(t_{k},x_{k}) \Delta
            + A_{t_{k+1}}^-G(t_{k},x_{k}) \Delta W_k\big)
            = 0.
      \end{align}
      On the one hand, we use $PA_t^- = A_t^-A_tA_t^- = A_t^-$ for any $t \in [0,T]$ to obtain $Pf(t,x) = f(t,x), x \in \mathbb{R}^{d}, t \in [0,T]$, and accordingly
      \begin{align*}
            Pu_{k+1}
            =&~
            \theta Pf(t_{k+1},u_{k+1}) \Delta
            + P A_{t_{k+1}} A_{t_k}^- x_k
            + (1-\theta) P A_{t_{k+1}}^-
              F(t_{k},x_{k}) \Delta \notag
            + PA_{t_{k+1}}^-G(t_{k},x_{k}) \Delta W_k
            \\=&~
            \theta f(t_{k+1},u_{k+1}) \Delta 
            + A_{t_{k+1}} A_{t_k}^- x_k
            + (1-\theta) A_{t_{k+1}}^-
              F(t_{k},x_{k}) \Delta \notag
            + A_{t_{k+1}}^-G(t_{k},x_{k}) \Delta W_k
            \\=&~
            u_{k+1} \in {\text{Im}}(P).
      \end{align*}
      On the other hand, 
      there exists a unique $v_{k+1} = \hat{V}(t_{k+1},u_{k+1}) \in \mathbb{R}^{d}$ such that 
      \begin{align*}
            A_{t_{k+1}}v_{k+1} = 0,
            \quad 
            RF(t_{k+1},u_{k+1}+v_{k+1})=0,
      \end{align*}
      which results in $Pv_{k+1} = A_{t_{k+1}}^-A_{t_{k+1}}v_{k+1} = 0$ and $Qv_{k+1} = (I-P)v_{k+1} = v_{k+1}$. Setting $x_{k+1} := u_{k+1} + v_{k+1}  \in \mathbb{R}^{d}$ allows us to get $Px_{k+1} = u_{k+1}$, $Qx_{k+1} = v_{k+1}$ and $x_{k+1} \in \mathcal{M}_{t_{k+1}}$. It follows from \eqref{lem:eq:u}, $RA_t = 0, RG(t,x) = 0, x  \in \mathbb{R}^{d}, t \in [0,T]$, $x_{k} \in \mathcal{M}_{t_{k}}$ and $A_{t_{k+1}}^- F(t_{k+1},x_{k+1}) = f(t_{k+1},u_{k+1})$ that
      \begin{align*}
            A_{t_{k+1}}^- H(x_{k+1})
            =&~
            A_{t_{k+1}}^-A_{t_{k+1}}x_{k+1} 
            - \theta A_{t_{k+1}}^-
            F(t_{k+1},x_{k+1}) \Delta 
            -
            \big( A_{t_{k+1}}^-A_{t_k}x_k
            \\&~+ (1-\theta) A_{t_{k+1}}^-
            F(t_{k},x_{k}) \Delta  
            + A_{t_{k+1}}^-G(t_k,x_k) 
            \Delta W_k\big)
            \\=&~
            u_{k+1} 
            - \theta f(t_{k+1},u_{k+1}) \Delta 
            -
            \big( A_{t_{k+1}}^-A_{t_k}x_k
            \\&~+ (1-\theta) A_{t_{k+1}}^-
            F(t_{k},x_{k}) \Delta  
            + A_{t_{k+1}}^-G(t_k,x_k) 
            \Delta W_k\big)
            \\=&~ 
            0,
      \end{align*}
      and 
      \begin{align*}
            R H(x_{k+1})
            =
            -\theta R F(t_{k+1},x_{k+1}) \Delta 
            -(1-\theta) R F(t_{k},x_{k}) \Delta
            =
            0.
      \end{align*}
      Along with the uniqueness of $u_{k+1}$ and $v_{k+1}$, $x_{k+1}$ is the unique solution to \eqref{lem:def:A^-H}, which means that there exists a unique solution $x_{k+1}  \in \mathcal{M}_{t_k}$ to \eqref{eq:ABEM}. Thus we complete the proof.\hfill$\square$

We continue to put the following assumption for SDAEs \eqref{eq:SDAE} adapting to the non-constant singular matrix $A_{t}$ for $t \in [0,T]$, which is not needed for the case of $A_{t} \equiv A,t \in [0,T]$ being a constant matrix; see Theorem \ref{th:conv2A}.

\begin{assumption}\label{asm:FGAt}
      There exist constants $L_2 > 0$ and $p_2 > 1$ such that for any $x,y \in \mathbb{R}^{d}$ and $t \in [0,T]$, 
      \begin{align}\label{asm:FGAt:rs}
            2\big\langle A_tx - A_ty, 
            F(t,x) - F(t,y) \big\rangle 
            +
            p_{2}\big|G(t,x) - G(t,y)\big|^{2}
            \leq 
            L_{2}\big|x - y\big|^{2}.
      \end{align}
\end{assumption}

\begin{rem}\label{rm:A1}
      While \eqref{asm:FG-FG} in Assumption \ref{asm:At} and \eqref{asm:FGAt:rs} appear similar, they are not equivalent. In fact, \eqref{asm:FGAt:rs} cannot be derived from \eqref{asm:FG-FG}. For example, take $d=3, r=2$ and consider
      \begin{equation*}
            A_t = {\rm{diag}}\bigg(\frac{1}{2},
            10, 0\bigg),
            \quad
            A_t^- = {\rm{diag}}\bigg(2, 
            \frac{1}{10}, 0\bigg),
            \quad
            P = {\rm{diag}}(1, 1, 0),
      \end{equation*}
      and
      \begin{equation*}
            F(t,x) = (-x_1^3, 0, 0)^{\top},
            \quad
            G(t,x) = {\rm{diag}}(0, x_1^2, 0),
            \quad x = (x_{1},x_{2},x_{3})^{\top} 
            \in \mathbb{R}^{3}.
      \end{equation*}
      Then for any $x,y \in \R^{3}$ and $t \in [0,T]$,
      \begin{align*}
            &~2\big\langle Px -Py, A^-_t F(t,x) 
            - A^-_t F(t,y) \big\rangle 
            + 10\big|A^-_t G(t,x) - A^-_t G(t,y)\big|^{2}
            \\=&~
            -4\big(x_1-y_1\big)\big(x_1^3-y_1^3\big) 
            +
            10\Big|\frac{1}{10}(x_1^2-y_1^2)\Big|^{2}
            \\=&~
            \big(x_1-y_1\big)^{2} 
            \Big(-4\big(x_1^2 + x_1y_1 + y_1^2\big)
            + \frac{1}{10}\big(x_1+y_1\big)^{2}\Big)
            \\\leq&~ 
            0,
      \end{align*}
      which implies that \eqref{asm:FG-FG} holds with $p_1=11$. Besides, one calculates that for any $x,y \in \R^{3}$ and $t \in [0,T]$, 
      \begin{align*}\label{rm:A1:AxF}
            &~2\big\langle A_tx - A_ty, 
            F(t,x) - F(t,y) \big\rangle 
            + \big|G(t,x) - G(t,y)\big|^{2}
            \\=&~
            -\big(x_1-y_1\big)
            \big(x_1^3-y_1^3\big) 
            +
            \big|x_1^2-y_1^2\big|^{2}
            \\=&~
            x_1 y_1 \big(x_1-y_1\big)^{2}.
      \end{align*}
      As there does not exist a constant $L_{2} > 0$ such that $x_1 y_1 \big(x_1-y_1\big)^{2} \leq L_{2}|x-y|^{2}$ for any $x,y \in \R^{3}$, Assumption \ref{asm:FGAt} thus fails to hold.
\end{rem}

The following lemma provides an upper mean square error bound via the
local truncated error, which is a key tool to establish the mean square convergence rate of STMs.

\begin{lemma}\label{lem:error}
      Suppose that Assumptions \ref{asm:At}, \ref{asm:FG} and \ref{asm:FGAt} hold and let $\theta \in [\frac12,1]$. If the stepsize $\Delta < \min\big\{
      \frac{1}{L_{1}\theta(1 + \hat{L}^{2})}, 
      \frac{1}{2L_{2}\theta}\big\}$, then there exists a constant $C > 0$, independent of $\Delta$, such that
      \begin{equation*}
            \max_{0 \leq k \leq K} 
            \mathbb{E}\big[\big|X_{t_k} 
            - x_k \big|^{2}\big]
            \leq 
            C\bigg(\sum_{k=1}^{K} \mathbb{E}
            \big[\big|\mathcal{R}_{k}\big|^{2}\big]
            +
            \Delta^{-1}\sum_{k=1}^{K} \mathbb{E}\big[ \big|\mathbb{E}\big( \mathcal{R}_k \,|\, \mathcal{F}_{t_{k-1}} \big)\big|^{2}
            \big]\bigg),
      \end{equation*}
      where 
      \begin{align*}
            \mathcal{R}_k
            := &~
            \theta \int_{t_{k-1}}^{t_k} 
            F(s,X_s) - F(t_k,X_{t_k}) {\rm{d}}s
            \\&~+
            (1-\theta) \int_{t_{k-1}}^{t_k} 
            F(s,X_s) - F(t_{k-1},X_{t_{k-1}}) {\rm{d}}s
            \\&~+ 
            \int_{t_{k-1}}^{t_k} G(s,X_s) 
            - G(t_{k-1},X_{t_{k-1}}) {\rm{d}}W_s,
            \quad k = 0,1,\cdots,K.
      \end{align*}
\end{lemma}

\textbf{Proof.}\quad
      For any $k = 0,1,\cdots,K$, letting $e_k := A_{t_k}X_{t_k} -A_{t_k}x_k$ and using \eqref{eq:SDAE} and \eqref{eq:ABEM} yield
      \begin{align*}
            e_k - \theta\big(F(t_{k},X_{t_{k}}) 
            - F(t_{k},x_{k}) \big)\Delta 
            =&~ 
            e_{k-1} 
            + (1-\theta) \big( F(t_{k-1},X_{t_{k-1}}) 
            - F(t_{k-1},x_{k-1}) \big)\Delta 
            \\&~
            + \big( G(t_{k-1},X_{t_{k-1}}) 
            - G(t_{k-1},x_{k-1}) \big) \Delta W_{k-1}
            + \mathcal{R}_{k}.
      \end{align*}
      Noting that $e_{k-1}$, $F(t_{k-1},X_{t_{k-1}})$ and $F(t_{k-1},x_{k-1})$, $G(t_{k-1},X_{t_{k-1}})$ and $G(t_{k-1},x_{k-1})$ are $\mathcal{F}_{t_{k-1}}$-measurable, and that $\Delta W_{k-1}$ is independent of $\mathcal{F}_{t_{k-1}}$, we use 
      \begin{align*}
            \mathbb{E}\big[\big\langle 
            e_{k-1} 
            +&~
            (1-\theta) \big(F(t_{k-1},X_{t_{k-1}}) 
            - F(t_{k-1},x_{k-1}) \big) \Delta,
            \\&~\big(G(t_{k-1},X_{t_{k-1}}) 
            - G(t_{k-1},x_{k-1})\big)\Delta W_{k-1}
            \big\rangle\big] = 0
      \end{align*}
      to get
      \begin{align}\label{lem:error:1clear}
            &~\mathbb{E}\big[\big| e_k 
            - \theta\big(F(t_{k},X_{t_{k}})    \notag
            - F(t_{k},x_{k})\big)\Delta \big|^{2}\big]
            \\=&~
            \mathbb{E}\big[\big| e_{k-1} 
            +
            (1-\theta) \big(F(t_{k-1},X_{t_{k-1}}) 
            - F(t_{k-1},x_{k-1}) \big) 
            \Delta \big|^{2}\big] \notag
            \\&~+
            \Delta \mathbb{E}\big[\big| 
            G(t_{k-1},X_{t_{k-1}}) 
            - G(t_{k-1},x_{k-1})\big|^{2}\big]
            +
            \mathbb{E}\big[\big| 
            \mathcal{R}_{k} \big|^{2}\big] \notag
            \\&~+
            2\mathbb{E}\big[\big\langle e_{k-1} 
            +  (1-\theta)\big(F(t_{k-1},X_{t_{k-1}}) 
            - F(t_{k-1},x_{k-1}) \big)\Delta, 
            \mathcal{R}_{k} \big\rangle\big] \notag
            \\&~+
            2\mathbb{E}\big[\big\langle 
            \big(G(t_{k-1},X_{t_{k-1}}) 
            - G(t_{k-1},x_{k-1})\big)\Delta W_{k-1},
            \mathcal{R}_{k} \big\rangle \big] \notag
        \\=&~
            \mathbb{E}\big[\big| e_{k-1} 
            - \theta \big(F(t_{k-1},X_{t_{k-1}}) 
            - F(t_{k-1},x_{k-1}) \big) \Delta 
            \big|^{2}\big] \notag
            \\&~+ 
            2\Delta \mathbb{E}\big[\big\langle e_{k-1}, F(t_{k-1},X_{t_{k-1}}) - F(t_{k-1},x_{k-1}) \big\rangle\big] \notag
            \\&~+ 
            (1-2\theta) \Delta^{2} \mathbb{E}
            \big[\big| F(t_{k-1},X_{t_{k-1}}) 
            - F(t_{k-1},x_{k-1}) \big|^{2}\big] \notag
            \\&~+
            \Delta \mathbb{E}\big[\big| 
            G(t_{k-1},X_{t_{k-1}}) 
            - G(t_{k-1},x_{k-1})\big|^{2}\big]
            +
            \mathbb{E}\big[\big| 
            \mathcal{R}_{k} \big|^{2}\big]  \notag
            \\&~+
            2\mathbb{E}\big[\big\langle e_{k-1} 
            +  (1-\theta)\big(F(t_{k-1},X_{t_{k-1}}) 
            - F(t_{k-1},x_{k-1}) \big)\Delta, 
            \mathcal{R}_{k} \big\rangle\big] \notag
            \\&~+
            2\mathbb{E}\big[\big\langle 
            \big(G(t_{k-1},X_{t_{k-1}}) 
            - G(t_{k-1},x_{k-1})\big)\Delta W_{k-1},
            \mathcal{R}_{k} \big\rangle \big].
      \end{align}
      As the random variables $e_{k-1}$, $F(t_{k-1},X_{t_{k-1}})$ and $F(t_{k-1},x_{k-1})$ are  $\mathcal{F}_{t_{k-1}}$-measurable, the Young inequlity and $\frac{1-\theta}{\theta}\leq 1$ indicate
      \begin{align}\label{lem:error:efRk}
            &~2\mathbb{E}\big[\big\langle e_{k-1} 
            + (1-\theta)\big(F(t_{k-1},X_{t_{k-1}}) 
            - F(t_{k-1},x_{k-1}) \big)\Delta, 
            \mathcal{R}_{k} \big\rangle\big]  \notag
            \\=&~  \notag
            2\mathbb{E}\big[\big\langle e_{k-1} 
            + (1-\theta)\big(F(t_{k-1},X_{t_{k-1}}) 
            - F(t_{k-1},x_{k-1}) \big)\Delta, 
            \mathbb{E}\big(\mathcal{R}_{k} \,|\, 
            \mathcal{F}_{t_{k-1}}\big\rangle\big]
            \\=&~  \notag
            2\frac{\theta-1}{\theta}
            \mathbb{E}\big[\big\langle e_{k-1} 
            - \theta\big(F(t_{k-1},X_{t_{k-1}}) 
            - F(t_{k-1},x_{k-1})\big)\Delta, 
            \\&~  \notag
            \mathbb{E}\big(\mathcal{R}_{k}
            \,|\, \mathcal{F}_{t_{k-1}} \big)
            \big\rangle\big]
            + 
            \frac{2}{\theta}\mathbb{E}\big[ 
            \big\langle e_{k-1}, 
            \mathbb{E}\big(\mathcal{R}_{k}
            \,|\, \mathcal{F}_{t_{k-1}} \big)
            \big\rangle\big]
            \\\leq&~  \notag
            \frac{1-\theta}{\theta}  \Big( 
            \Delta \mathbb{E}\big[|e_{k-1} 
            - \theta \big(F(t_{k-1},X_{t_{k-1}}) 
            - F(t_{k-1},x_{k-1}) \big) 
            \Delta\big|^{2}\big]
            \\&~+  \notag
            \Delta^{-1} \mathbb{E}\big[\big|
            \mathbb{E}\big(\mathcal{R}_{k}
            \,|\, \mathcal{F}_{t_{k-1}} \big) 
            \big|^{2}\big]\Big)
            +
            \Delta \mathbb{E}\big[
            \big|e_{k-1}\big|^{2}\big]
            + 
            \frac{1}{\theta^2} \Delta^{-1} \mathbb{E}
            \big[\big|\mathbb{E}\big(\mathcal{R}_{k}\,|\, 
            \mathcal{F}_{t_{k-1}} \big) \big|^{2}\big] 
        \\\leq&~
            \Delta \mathbb{E}\big[|e_{k-1} 
            - \theta \big(F(t_{k-1},X_{t_{k-1}}) 
            - F(t_{k-1},x_{k-1}) \big) 
            \Delta\big|^{2}\big]
            \\&~+
            \Delta \mathbb{E}\big[
            \big|e_{k-1}\big|^{2}\big] \notag
            + 
            \frac{\theta^{2}+1}{\theta^{2}} 
            \Delta^{-1} \mathbb{E}\big[\big|
            \mathbb{E}\big(\mathcal{R}_{k} \,|\, 
            \mathcal{F}_{t_{k-1}} \big) \big|^{2}\big].
      \end{align}
      By means of the Young inequlity again, we have
      \begin{align}\label{lem:error:leq0}
            &~2\mathbb{E}\big[\big\langle 
            \big(G(t_{k-1},X_{t_{k-1}}) 
            - G(t_{k-1},x_{k-1})\big) \Delta W_{k-1}, \mathcal{R}_k \big\rangle \big] \notag
            \\\leq&~
            \big(p_{2}-1\big) \Delta 
            \mathbb{E}\big[\big|G(t_{k-1},X_{t_{k-1}}) 
            - G(t_{k-1},x_{k-1}) \big|^{2}\big]
            +
            \frac{1}{p_2-1} \mathbb{E}\big[
            \big|\mathcal{R}_k \big|^{2}\big].
      \end{align}
      Inserting \eqref{lem:error:efRk} and \eqref{lem:error:leq0} into \eqref{lem:error:1clear}
      and utilizing Assumption \ref{asm:FGAt} as well as  $1-2\theta \leq 0$ lead to
      \begin{align}\label{lem:error:1clear1}
            &~\mathbb{E}\big[\big| e_k 
            - \theta\big(F(t_{k},X_{t_{k}}) \notag
            - F(t_{k},x_{k})\big)\Delta \big|^{2}\big]
        \\\leq&~
            \mathbb{E}\big[\big| e_{k-1} 
            - \theta \big(F(t_{k-1},X_{t_{k-1}}) 
            - F(t_{k-1},x_{k-1}) \big) \Delta 
            \big|^{2}\big] \notag
            \\&~+ 
            L_{2}\Delta \mathbb{E}\big[\big|
            X_{t_{k-1}} - x_{k-1}) \big|^{2}\big] \notag
            -
            p_{2}\Delta \mathbb{E}\big[\big| 
            G(t_{k-1},X_{t_{k-1}}) 
            - G(t_{k-1},x_{k-1}) \big|^{2}\big] \notag
            \\&~+
            \Delta \mathbb{E}\big[\big| 
            G(t_{k-1},X_{t_{k-1}}) 
            - G(t_{k-1},x_{k-1})\big|^{2}\big]
            +
            \mathbb{E}\big[\big| 
            \mathcal{R}_{k} \big|^{2}\big]  \notag
            \\&~+
            \Delta \mathbb{E}\big[|e_{k-1} 
            - \theta \big(F(t_{k-1},X_{t_{k-1}}) 
            - F(t_{k-1},x_{k-1}) \big) 
            \Delta\big|^{2}\big] \notag
            \\&~+
            \Delta \mathbb{E}\big[
            \big|e_{k-1}\big|^{2}\big] \notag
            + 
            \frac{\theta^{2}+1}{\theta^{2}} 
            \Delta^{-1} \mathbb{E}\big[\big|
            \mathbb{E}\big(\mathcal{R}_{k} \,|\, 
            \mathcal{F}_{t_{k-1}} \big) \big|^{2}\big]
            \\&~+
            \big(p_{2}-1\big) \Delta 
            \mathbb{E}\big[\big|G(t_{k-1},X_{t_{k-1}}) 
            - G(t_{k-1},x_{k-1}) \big|^{2}\big]
            +
            \frac{1}{p_2-1} \mathbb{E}\big[
            \big|\mathcal{R}_k \big|^{2}\big] \notag
        \\=&~
            \big(1 + \Delta\big)
            \mathbb{E}\big[\big| e_{k-1} 
            - \theta \big(F(t_{k-1},X_{t_{k-1}}) 
            - F(t_{k-1},x_{k-1}) \big) \Delta 
            \big|^{2}\big]
            +
            \Delta \mathbb{E}\big[
            \big|e_{k-1}\big|^{2}\big] \notag
            \\&~+
            L_{2}\Delta \mathbb{E}\big[\big|
            X_{t_{k-1}} - x_{k-1}\big|^{2}\big]
            +
            \frac{p_2}{p_2-1} \mathbb{E}\big[
            \big|\mathcal{R}_k \big|^{2}\big]
            + 
            \frac{\theta^{2}+1}{\theta^{2}} 
            \Delta^{-1} \mathbb{E}\big[\big|
            \mathbb{E}\big(\mathcal{R}_{k} \,|\, 
            \mathcal{F}_{t_{k-1}} \big) \big|^{2}\big].
      \end{align}
      Lemma \ref{lem:def} ensures $x_{k} \in \mathcal{M}_{t_{k}}$ for any $k = 0,1,2, \cdots, K$, which means $x_{k} = Px_{k} + \hat{V}(t_{k},Px_{k}) = u_{k} + \hat{V}(t_{k},u_{k})$. Along with \eqref{lem-v-rs3}, $P = A_{t_k}^{-}A_{t_k}$ and Lemma \ref{lem:prop:A}, one gets that for any $k = 0,1,\cdots,K$,
      \begin{align}\label{lem:error:final}
            \mathbb{E}\big[\big|X_{t_k} 
            - x_{k} \big|^{2} \big]
            =&~ \notag
            \mathbb{E}\big[\big|\big(PX_{t_k} +
            \hat{V}(t_k,PX_{t_k})\big)
            - \big(Px_k + \hat{V}(t_k,Px_k)\big)  
            \big|^{2}\big]
            \\\leq&~ \notag
            2\big(1 + \hat{L}\big)^{2}
            \mathbb{E}\big[\big|PX_{t_k} 
            - Px_k \big|^{2}\big]
            =
            2\big(1 + \hat{L}\big)^{2}
            \mathbb{E}\big[\big|A_{t_k}^-
            e_k \big|^{2}\big]
            \\=&~
            2\big(1 + \hat{L}\big)^{2}
            \big|A_{t_k}^-\big|^{2}
            \mathbb{E}\big[\big|
            e_k \big|^{2}\big]
            \leq 
            2r\Big(\frac{1 + \hat{L}}
            {\underline{\sigma}}\Big)^{2}
            \mathbb{E}\big[\big|
            e_k \big|^{2}\big],
      \end{align}
      which in combination with \eqref{lem:error:1clear1} yields 
      \begin{align*}
            &~\mathbb{E}\big[\big| e_k 
            - \theta\big(F(t_{k},X_{t_{k}}) \notag
            - F(t_{k},x_{k})\big)\Delta \big|^{2}\big]
        \\\leq&~ \notag
            \big(1 + \Delta\big)
            \mathbb{E}\big[\big| e_{k-1} 
            - \theta \big(F(t_{k-1},X_{t_{k-1}}) 
            - F(t_{k-1},x_{k-1}) \big) \Delta 
            \big|^{2}\big]
            \\&~+
            \bigg(1 + 2rL_{2}\Big(\frac{1 + \hat{L}}
            {\underline{\sigma}}\Big)^{2}\bigg)
            \Delta \mathbb{E}\big[
            \big|e_{k-1}\big|^{2}\big] \notag
            +
            \frac{p_2}{p_2-1} \mathbb{E}\big[
            \big|\mathcal{R}_k \big|^{2}\big]
            \\&~+ 
            \frac{\theta^{2}+1}{\theta^{2}} 
            \Delta^{-1} \mathbb{E}\big[\big|
            \mathbb{E}\big(\mathcal{R}_{k} \,|\, 
            \mathcal{F}_{t_{k-1}} \big) \big|^{2}\big].
      \end{align*}
      Noting that $e_{0} = 0$, $F(t_{0},X_{t_{0}}) - F(t_{0},x_{0}) = 0$ and $(1+\Delta)^{k} \leq e^{T}$ for $k = 0,1,\cdots,K$, an iterative argument gives 
      \begin{align*}
            &~\mathbb{E}\big[\big| e_k 
            - \theta\big(F(t_{k},X_{t_{k}}) \notag
            - F(t_{k},x_{k})\big)\Delta \big|^{2}\big]
        \\\leq&~ \notag
            \big(1 + \Delta\big)^{k} 
            \mathbb{E}\big[\big| e_{0} 
            - \theta\big(F(t_{0},X_{t_{0}})
            - F(t_{0},x_{0})\big)\Delta \big|^{2}\big]
            \\&~+ \notag
            \bigg(1 + 2rL_{2}\Big(\frac{1 + \hat{L}}
            {\underline{\sigma}}\Big)^{2}\bigg)
            \Delta \sum_{i=0}^{k-1} (1+\Delta)^{k-1-i} \mathbb{E}\big[\big|e_{i}\big|^{2}\big]
            \\&~+ \notag
            \frac{p_2}{p_2-1} \sum_{i=1}^{k} 
            (1+\Delta)^{k-i} \mathbb{E}\big[
            \big|\mathcal{R}_{i} \big|^{2}\big]
            +
            \frac{\theta^{2}+1}{\theta^{2}} 
            \Delta^{-1} \sum_{i=1}^{k} 
            (1+\Delta)^{k-i} \mathbb{E}\big[\big|
            \mathbb{E}\big(\mathcal{R}_{i} \,|\, 
            \mathcal{F}_{t_{i-1}} \big) \big|^{2}\big]
        \\\leq&~ \notag
            \bigg(1 + 2rL_{2}\Big(\frac{1 + \hat{L}}
            {\underline{\sigma}}\Big)^{2}\bigg) e^{T}
            \Delta \sum_{i=0}^{k-1}  \mathbb{E}\big[\big|e_{i}\big|^{2}\big]
            +
            \frac{p_{2} e^{T}}{p_{2}-1} 
            \sum_{i=1}^{k} \mathbb{E}\big[
            \big|\mathcal{R}_{i} \big|^{2}\big]
            \\&~+
            \frac{\theta^{2}+1}{\theta^{2}} 
            e^{T} \Delta^{-1} \sum_{i=1}^{k} 
            \mathbb{E}\big[\big|
            \mathbb{E}\big(\mathcal{R}_{i} \,|\, 
            \mathcal{F}_{t_{i-1}} \big) \big|^{2}\big].
      \end{align*}
      Assumption \ref{asm:FGAt} indicates
      \begin{align*}
            \big|e_{k} - \theta\big(F(t_{k},X_{t_k}) 
            - F(t_{k},x_{k}) \big)\Delta\big|^{2}
            \geq
            \big|e_{k}\big|^{2} 
            - 
            2\theta\Delta \big\langle e_k,
            F(t_{k},X_{t_k}) 
            - F(t_{k},x_{k}) \big\rangle
            \geq
            \big(1 - L_{2}\theta\Delta\big)|e_{k}|^{2},
      \end{align*}
      and consequently
      \begin{align*}
            \big(1 - L_{2}\theta\Delta) 
            \mathbb{E}\big[\big|e_{k}\big|^{2}\big]
            \leq&~ \notag
            \bigg(1 + 2rL_{2}\Big(\frac{1 + \hat{L}}
            {\underline{\sigma}}\Big)^{2}\bigg) e^{T}
            \Delta \sum_{i=0}^{k-1}  \mathbb{E}\big[\big|e_{i}\big|^{2}\big]
            \\&~+
            \frac{p_{2} e^{T}}{p_{2}-1} 
            \sum_{i=1}^{K} \mathbb{E}\big[
            \big|\mathcal{R}_{i} \big|^{2}\big]
            +
            \frac{\theta^{2}+1}{\theta^{2}} 
            e^{T} \Delta^{-1} \sum_{i=1}^{K} 
            \mathbb{E}\big[\big|
            \mathbb{E}\big(\mathcal{R}_{i} \,|\, 
            \mathcal{F}_{t_{i-1}} \big) \big|^{2}\big].
      \end{align*}
      Making using of $\Delta<\frac{1}{2L_{2}\theta}$ and the discrete Gronwall inequality, there exists a constant $C > 0$, independent of $\Delta$, such that for any $k = 0,1,\cdots,K$, 
      \begin{align*}
            \mathbb{E}\big[\big|e_{k}\big|^{2}\big]
            \leq
            C\bigg(\sum_{i=1}^{K} \mathbb{E}\big[
            \big|\mathcal{R}_{i} \big|^{2}\big]
            +
            \Delta^{-1} \sum_{i=1}^{k} 
            \mathbb{E}\big[\big|
            \mathbb{E}\big(\mathcal{R}_{i} \,|\, 
            \mathcal{F}_{t_{i-1}} \big)
            \big|^{2}\big]\bigg),
      \end{align*}
      which in combination with \eqref{lem:error:final} gives the desired result.\hfill$\square$

With the help of the previously established upper mean square error
bound, we are able to reveal the mean square convergence rates of STMs, which is one of the main results of this work.

\begin{theorem}\label{thm:convAt}
      Suppose that Assumptions \ref{asm:At}, \ref{asm:FG} and \ref{asm:FGAt} hold and let $\theta \in [\frac12,1]$. If the stepsize $\Delta < \min\big\{\frac{1}{L_{1}\theta(1 + \hat{L}^{2})}, \frac{1}{2L_{2}\theta}\big\}$, then there exists a constant $C > 0$, independent of $\Delta$, such that
      \begin{equation*}
            \max_{0\leq k \leq K} 
            \mathbb{E}\big[\big|X_{t_k} 
            - x_{k} \big|^{2}\big]
            \leq C \Delta.
      \end{equation*}
\end{theorem}

\textbf{Proof.}\quad
      By \eqref{asm:FG-absF} and \eqref{asm:FG-absG}, we have that for any $s,t \in [0,T]$,
      \begin{align*}
            &~\big|F(t,X_t) - F(s,X_s)\big|^{2} 
            + \big|G(t,X_t) - G(s,X_s)\big|^{2}
            \\\leq&~
            C\Big(\big(1 + |X_t| + |X_s|\big)^{\gamma-1} 
            \big|X_t-X_s\big| 
            + \big(1 + |X_t| + |X_s|\big)^{\gamma} 
            \big|t-s\big|\Big)^{2}
            \\&~+ 
            C\Big(\big(1 + |X_t| + |X_s|\big)^{\gamma-1} \big|X_t-X_s\big|^{2} 
            + \big(1 + |X_t| + |X_s|\big)^{\gamma+1} 
            \big|t-s\big|^{2}\Big)
            \\\leq&~
            C\Big(\big(1 + |X_t| + |X_s|\big)^{2\gamma-2}
            \big|X_t-X_s\big|^{2} 
            + \big(1 + |X_t| + |X_s|\big)^{2\gamma} 
            \big|t-s\big|^{2}\Big).
      \end{align*}
      It follows from the H\"older inequlity, $2\gamma\leq 4\gamma-2 <p_1$, and Lemmas \ref{lem:bound}, \ref{lem:UX} that
      \begin{align*}
            &~\mathbb{E}\big[\big|F(t,X_t)
            - F(s,X_s) \big|^{2}\big] 
            + \mathbb{E}\big[\big|G(t,X_t) 
            - G(s,X_s) \big|^{2}\big]
            \\\leq&~
            C\Big(1 
            + 
            \mathbb{E}\big[\big|X_t|^{4\gamma-2}\big] 
            +
            \mathbb{E}\big[\big|X_s|^{4\gamma-2}\big]
            \Big)^{\frac{\gamma-1}{2\gamma-1}}
            \Big(\mathbb{E}\big[\big|X_t 
            - X_s\big|^{\frac{4\gamma-2}{\gamma}}\big]
            \Big)^{\frac{\gamma}{2\gamma-1}}
            \\&~+ 
            C|t-s|^{2}\Big(1 
            + \mathbb{E}\big[\big|X_t|^{2\gamma}\big]
            + \mathbb{E}\big[\big|X_s|^{2\gamma}\big]\Big)
            \\\leq&~ 
            C|t-s|.
      \end{align*}
      By the H\"{o}lder inequlity, $\theta \in [1/2,1]$ and the It\^{o} isometry, one deduces that for any $k = 1,2, \cdots, K$,
      \begin{align}\label{lem:conv:r1}
            &~\mathbb{E}\big[\big|\mathcal{R}_{k}
            \big|^{2}\big] \notag
            \\\leq&~ \notag
            3\mathbb{E}\bigg[\Big|
            \int_{t_{k-1}}^{t_k} F(s,X_s) 
            - F(t_k,X_{t_k}) \text{d}s 
            \Big|^{2}\bigg]
            +
            3\mathbb{E}\bigg[\Big|
            \int_{t_{k-1}}^{t_k} F(s,X_s) 
            - F(t_{k-1},X_{t_{k-1}}) 
            \text{d}s \Big|^{2}\bigg]
            \\&~+ \notag
            3\mathbb{E}\bigg[\Big|
            \int_{t_{k-1}}^{t_k} G(s,X_s) 
            - G(t_{k-1},X_{t_{k-1}}) 
            \text{d}W_s \Big|^{2}\bigg]
            \\\leq&~ \notag
            3\Delta \int_{t_{k-1}}^{t_k}
            \mathbb{E}\big[\big|F(s,X_s) 
            - F(t_k,X_{t_k}) \big|^{2}\big] \text{d}s 
            +
            3\Delta \int_{t_{k-1}}^{t_k} 
            \mathbb{E}\big[\big|F(s,X_s) 
            - F(t_{k-1},X_{t_{k-1}}) 
            \big|^{2}\big] \text{d}s
            \\&~ \notag
            +3 \int_{t_{k-1}}^{t_k} 
            \mathbb{E}\big[\big|G(s,X_s) 
            - G(t_{k-1},X_{t_{k-1}}) 
            \big|^{2}\big] \text{d}s
            \\
            \leq &~ C\Delta^2,
      \end{align}
      and 
      \begin{align}\label{lem:conv:r2}
            &~\mathbb{E}\big[\big|\mathbb{E}
            \big( \mathcal{R}_{k} \,|\, 
            \mathcal{F}_{t_{k-1}} \big) \big|^{2}\big] \notag
            \\=&~   \notag
            \mathbb{E}\bigg[\Big|
            \theta \int_{t_{k-1}}^{t_k} 
            F(s,X_s) - F(t_k,X_{t_k}) \text{d}s
            +(1-\theta) \int_{t_{k-1}}^{t_k} 
            F(s,X_s) - F(t_{k-1},X_{t_{k-1}}) 
            \text{d}s \Big|^{2}\bigg]
            \\\leq&~   \notag
            2\Delta \int_{t_{k-1}}^{t_k} 
            \mathbb{E}\big[\big|F(s,X_s) 
            - F(t_k,X_{t_k}) \big|^{2}\big]\text{d}s 
            +
            2\Delta \int_{t_{k-1}}^{t_k} 
            \mathbb{E}\big[\big|F(s,X_s) 
            - F(t_{k-1},X_{t_{k-1}}) 
            \big|^{2}\big]\text{d}s  \notag
            \\\leq&~  
            C\Delta^{3}.
      \end{align}
      Together with Lemma \ref{lem:error}, we obtain the required result and finish the proof.\hfill$\square$

Generally speaking, \eqref{asm:FGAt:rs} cannot be derived from \eqref{asm:FG-FG}; see Remark \ref{rm:A1}. However, this result holds true in a special case, as demonstrated below.

\begin{prop}\label{prop:At}
      If Assumption \ref{asm:At} holds with $r = 1$, then \eqref{asm:FG-FG} implies \eqref{asm:FGAt:rs} with $p_{2} = p_{1} - 1$ and $L_{2} = \bar{\sigma}^{2}L_{1}$.
\end{prop}

\textbf{Proof.}\quad
      Since 
      $A_t = M \mathrm{diag}(\sigma_{1}(t), 0, \cdots, 0) N$ with $M$ and $N$ being $d$-order orthogonal matrices, we have 
      \begin{align*}
            A_{t}^{-} 
            =
            N^{\top} \mathrm{diag}\bigg(
            \frac{1}{\sigma_{1}(t)}, 0, 
            \cdots, 0\bigg) M^{\top} 
            =
            \frac{1}{\sigma_{1}^{2}(t)} A_{t}^{\top},
            \quad t \in [0,T].
      \end{align*}
      For any $x,y \in \mathbb{R}^{d}$ and $t \in [0,T]$, utilizing $P A_t^- = \big(A_t^-A_t\big) A_t^- = A_t^-$ leads to
      \begin{align*}
            \big\langle A_tx - A_ty, F(t,x) -F(t,y) \big\rangle
            =&~
            \big\langle x - y, A_t^{\top}
            \big(F(t,x) - F(t,y)\big) \big\rangle
            \\=&~
            \sigma_{1}^2(t) \big\langle x - y, 
            A_t^- (F(t,x) -F(t,y)) \big\rangle
            \\=&~
            \sigma_{1}^2(t) \big\langle x - y, P A_t^- 
            \big(F(t,x) - F(t,y)\big) \big\rangle
            \\=&~
            \sigma_{1}^2(t) \big\langle P(x-y), 
            A_t^- \big(F(t,x) - F(t,y)\big) \big\rangle.
      \end{align*}
      We then combine with $R G(t,x) = 0$, $R = I - A_tA_t^-$ and and $A_{t}^{\top} = \sigma_{1}^{2}(t)A_{t}^{-}$ to get
      \begin{align*}
            \big| G(t,x) - G(t,y) \big|^{2} 
            =&~
            \big| \big(G(t,x) - G(t,y)\big)^{\top} 
                I \big(G(t,x) - G(t,y)\big) \big|
            \\=&~
            \big| \big(G(t,x) -G(t,y)\big)^{\top} 
            \big(R + A_tA_t^-\big) 
            \big(G(t,x) - G(t,y)\big) \big|
            \\=&~
            \big|\big(G(t,x) - G(t,y)\big)^{\top} 
            A_tA_t^- \big(G(t,x) - G(t,y)\big) \big|
            \\=&~
            \sigma_{1}^{2}(t) \big| \big(G(t,x) - G(t,y)\big)^{\top} 
            (A_t^-)^{\top} A_t^- \big(G(t,x) - G(t,y)\big) \big|
            \\=&~
            \sigma_{1}^2(t) 
            \big| A_t^-\big(G(t,x) - G(t,y)\big) \big|^{2}.
      \end{align*}
      It follows from \eqref{asm:FG-FG} and Assumption \ref{asm:At} that
      \begin{align*}
            &~2\big\langle A_tx - A_ty, 
            F(t,x) - F(t,y) \big\rangle 
            + 
            \big(p_{1} - 1\big)
            \big|G(t,x) - G(t,y)\big|^{2}
            \\=&~
            \sigma_{1}^2(t) \Big(2\big\langle Px -Py, 
            A^-_t F(t,x) - A^-_t F(t,y) \big\rangle 
            +
            \big(p_{1} - 1\big) \big|A^-_t G(t,x) 
            - A^-_t G(t,y)\big|^{2}\Big)
            \\=&~
            \sigma_{1}^{2}(t) L_{1}\big|x-y\big|^{2}
            \leq
            \bar{\sigma}^{2} L_{1}\big|x-y\big|^{2},
      \end{align*}
      which implies \eqref{asm:FGAt:rs} holds with $p_{2} = p_{1} - 1$ and $L_{2} = L_{1} \bar{\sigma}^{2}$. Thus we complete the proof. \hfill$\square$

As a direct consequence of Theorem \ref{thm:convAt} and Proposition \ref{prop:At}, we obtain the following corollary, whose proof is thus omitted.
\begin{corollary}\label{co:conv}
      Suppose that Assumption \ref{asm:At} holds with $r=1$ and Assumption \ref{asm:FG} holds. If $\theta \in [\frac12,1]$ and $\Delta< \frac{1}{L_{1}\theta(1 + \hat{L}^{2})}$, then there exists a constant $C > 0$, independent of $\Delta$, such that
      \begin{equation*}
            \max_{0 \leq k \leq K} 
            \mathbb{E}\big[\big|X_{t_k} 
            - x_k \big|^{2}\big]
            \leq 
            C\Delta.
      \end{equation*}
\end{corollary}

\section{Mean square convergence rates of STM with constant $A $}\label{sec:orderA}
Although Theorem \ref{thm:convAt} can directly gives the convergence rate for STMs of \eqref{eq:SDAE} with a constant matrix $A_{t} \equiv A,t \in [0,T]$,
Assumption \ref{asm:FGAt} is redundant. In fact, in this section we aim to achieve the same convergence rate without Assumption \ref{asm:FGAt}. To this end, we begin with an upper mean square error bound as done in Lemma \ref{lem:error}.

\begin{lemma}\label{lem:error2}
      Let $A_t \equiv A, t \in [0,T]$ be a constant matrix and suppose that Assumption \ref{asm:FG} holds. If $\theta \in [\frac12,1]$ and $\Delta < \frac{1}{2L_{1}\theta(1 + \hat{L}^{2})}$, then there exists a constant $C > 0$, independent of $\Delta$ such that
      \begin{equation*}
            \max_{0\leq k \leq K} 
            \mathbb{E}\big[\big|X_{t_k} 
            - x_k \big|^{2}\big]
            \leq 
            C\bigg(\sum_{k=1}^{K}\mathbb{E}\big[ 
            \big| \mathcal{Q}_{k} \big|^{2}\big]
            +
            \Delta^{-1} \sum_{k=1}^{K} 
            \mathbb{E}\big[ \mathbb{E} \big( 
            \mathcal{Q}_{k} \,|\, \mathcal{F}_{t_{k-1}} 
            \big) \big|^{2} \big]\bigg),
      \end{equation*}
      where 
      \begin{align*}
            \mathcal{Q}_k
            :=&~
            \theta \int_{t_{k-1}}^{t_k} f(s,U_s) 
            - f(t_k,U_{t_k}) {\rm{d}}s
            \\&~+
            (1-\theta) \int_{t_{k-1}}^{t_k} f(s,U_s) 
            - f(t_{k-1},U_{t_{k-1}}) {\rm{d}}s
            \\&~+ 
            \int_{t_{k-1}}^{t_k} g(s,U_s) 
            - g(t_{k-1},U_{t_{k-1}}) {\rm{d}}W_{s},
            \quad k = 1,2, \cdots, K.
      \end{align*}
\end{lemma}

\textbf{Proof.}\quad
      Since \eqref{def:Ff} implies $f(t,Px) = A^-F(t,x)$ and $g(t,Px)=A^-G(t,x)$ for any $x \in \mathcal{M}_{t}, t \in [0,T]$, we use
      $X_{t_k}, x_k \in \mathcal{M}_{t_k}$ and $U_{t_k} = PX_{t_k}, u_k = Px_k$ for $k = 0,1,\cdots,K$, and \eqref{eq:ABEM} to get
      \begin{align*}
            u_{k+1} 
            =
            u_k + \theta f(t_{k+1},u_{k+1}) \Delta 
            + (1-\theta) f(t_{k},u_{k}) \Delta 
            + g(t_k,u_k) \Delta {W_{k}}.
      \end{align*}
      For any $k = 0,1,\cdots,K$, letting $e_k := X_{t_k} -x_k$ and applying \eqref{eq:inherent_U} give $Pe_k := U_{t_k} -u_k$ and 
      \begin{align*}
            Pe_{k} - \theta \big(f(t_{k},U_{t_k}) 
            - f(t_{k},u_{k}) \big) \Delta 
            =&~ 
            Pe_{k-1} +
            (1-\theta) \big(f(t_{k-1},U_{t_{k-1}}) 
            - f(t_{k-1},u_{k-1}) \big) \Delta 
            \\&~+ 
            \big(g(t_{k-1},U_{t_{k-1}}) 
            - g(t_{k-1},u_{k-1})\big) \Delta W_{k-1}
            + \mathcal{Q}_{k}.
      \end{align*}
      Similar to \eqref{lem:error:1clear}, the $\mathcal{F}_{t_{k-1}}$-measurability of $e_{k-1}$, $U_{t_{k-1}}$ and $u_{k-1}$ shows
      \begin{align}\label{lem:error2:1clear}
            &~\mathbb{E}\big[\big|Pe_k  
            - \theta\big(f(t_{k},U_{t_{k}}) 
            - f(t_{k},u_{k})\big)\Delta\big|^{2}\big]\notag
        \\=&~ \notag
            \mathbb{E}\big[\big|Pe_{k-1} 
            - \theta\big(f(t_{k-1},U_{t_{k-1}}) 
            - f(t_{k-1},u_{k-1})\big)\Delta
            \big|^{2}\big]
            \\&~ \notag
            + 2\Delta \mathbb{E}\big[ 
            \big\langle Pe_{k-1}, f(t_{k-1},U_{t_{k-1}}) 
            - f(t_{k-1},u_{k-1}) \big\rangle\big]
            \\&~ \notag
            + (1-2\theta) \Delta^{2} 
            \mathbb{E}\big[\big|f(t_{k-1},U_{t_{k-1}}) 
            - f(t_{k-1},u_{k-1}) \big|^{2}\big]
            \\&~+ \notag
            \Delta \mathbb{E}\big[\big|
            g(t_{k-1},U_{t_{k-1}}) 
            - g(t_{k-1},u_{k-1})\big|^{2}\big]
            + 
            \mathbb{E}\big[\big|
            \mathcal{Q}_k \big|^{2}\big]
            \\&~+ \notag
            2\mathbb{E}\big[\big\langle 
            Pe_{k-1} + (1-\theta)\big(
            f(t_{k-1},U_{t_{k-1}}) 
            - f(t_{k-1},u_{k-1})\big)\Delta, 
            \mathcal{Q}_{k} \big\rangle\big]
            \\&~+
            2\mathbb{E}\big[\big\langle 
            \big(g(t_{k-1},U_{t_{k-1}}) 
            - g(t_{k-1},u_{k-1})\big)\Delta W_{k-1}, 
            \mathcal{Q}_{k} \big\rangle\big].
      \end{align}
      Utilizing \eqref{asm:FG-FG} and \eqref{def:Ff} ensures that for any $x,y \in \mathcal{M}_{t}, t \in [0,T]$,
      \begin{align}\label{asm:FG-FG2}
            2\big\langle Px-Py, 
            f(t,Px) - f(t,Py) \big\rangle 
            +
            (p_{1}-1)
            \big|g(t,Px) - g(t,Py))\big|^{2}
            \leq 
            L_{1}\big|x-y\big|^{2},
      \end{align}
      and accordingly
      \begin{align*}
             &~2\Delta \mathbb{E}\big[ 
            \big\langle Pe_{k-1}, f(t_{k-1},U_{t_{k-1}}) 
            - f(t_{k-1},u_{k-1}) \big\rangle\big]
            \\\leq&~
            L_{1}\Delta\mathbb{E}\big[
            \big|e_{k-1}\big|^{2}\big]
            -
            (p_{1}-1)\Delta\mathbb{E}\big[\big|
            \big|g(t_{k-1},U_{t_{k-1}}) 
            - g(t_{k-1},u_{k-1})\big|^{2}\big].
      \end{align*}
      Analogously to \eqref{lem:error:efRk} and \eqref{lem:error:leq0}, one can derive 
      \begin{align*}
            &~2\mathbb{E}\big[\big\langle 
            Pe_{k-1} + (1-\theta)\big(
            f(t_{k-1},U_{t_{k-1}}) 
            - f(t_{k-1},u_{k-1})\big)\Delta, 
            \mathcal{Q}_{k} \big\rangle\big]
            \\\leq&~             
            \Delta \mathbb{E}\big[\big|Pe_{k-1} 
            - \theta\big(f(t_{k-1},U_{t_{k-1}}) 
            - f(t_{k-1},u_{k-1})\big) \Delta
            \big|^{2}\big]
            \\&~+ 
            \Delta \mathbb{E}\big[\big|
            Pe_{k-1}\big|^{2}\big]
            + 
            \frac{\theta^{2} + 1}{\theta^{2}}
            \Delta^{-1} \mathbb{E}\big[\big| \mathbb{E} \big(\mathcal{Q}_{k} \,|\, 
            \mathcal{F}_{t_{k-1}}\big)\big|^{2}\big],
      \end{align*}
      and 
      \begin{align*}
            &~2\mathbb{E}\big[\big\langle 
            \big(g(t_{k-1},U_{t_{k-1}}) 
            - g(t_{k-1},u_{k-1})\big)\Delta W_{k-1}, \mathcal{Q}_k \big\rangle\big]
            \\\leq&~
            (p_{1}-2) \Delta \mathbb{E}\big[\big|
            g(t_{k-1},U_{t_{k-1}}) 
            - g(t_{k-1},u_{k-1})\big|^{2}\big]
            +
            \frac{1}{p_{1}-2} \mathbb{E}\big[\big|
            \mathcal{R}_{k} \big|^{2}\big].
      \end{align*}
      It follows from \eqref{lem:error2:1clear} and $1-2\theta \leq 0$ that   
      \begin{align*}
            &~\mathbb{E}\big[\big|Pe_k  
            - \theta\big(f(t_{k},U_{t_{k}}) 
            - f(t_{k},u_{k})\big)\Delta\big|^{2}\big]\notag
        \\\leq&~ \notag
            (1 + \Delta)
            \mathbb{E}\big[\big|Pe_{k-1} 
            - \theta\big(f(t_{k-1},U_{t_{k-1}}) 
            - f(t_{k-1},u_{k-1})\big)\Delta
            \big|^{2}\big]
            \\&~+ 
            \big(1 + 2L_{1}(1+\hat{L}^{2})\big)
            \Delta \mathbb{E}\big[
            \big|Pe_{k-1}\big|^{2}\big]
            +
            \frac{p_1-1}{p_1-2}\mathbb{E}\big[
            \big|\mathcal{Q}_{k}\big|^{2}\big] \notag
            \\&~+ 
            \frac{\theta^{2} + 1}{\theta^{2}}
            \Delta^{-1} \mathbb{E}\big[\big| 
            \mathbb{E} \big(\mathcal{Q}_{k} \,|\, 
            \mathcal{F}_{t_{k-1}}\big)\big|^{2}\big],
      \end{align*}
      where we have used the fact that for $k = 0,1,\cdots,K$,
      \begin{align}\label{lem:error2:final}
            \big|e_{k}\big|^{2}
            =
            \big|\big(U_{t_k} + \hat{V}(t_k,U_{t_k})\big) 
            - \big(u_k + \hat{V}(t_k,u_k)\big) \big|^{2}
            \leq
            2(1 + \hat{L}^{2}) \big|Pe_{k}\big|^{2}
      \end{align}
      due to \eqref{lem-v-uv}. By iteration, we obtain \begin{align}\label{lem:error2:leq3}
            &~\mathbb{E}\big[\big|Pe_k  
            - \theta\big(f(t_{k},U_{t_{k}}) 
            - f(t_{k},u_{k})\big)\Delta\big|^{2}\big]\notag
            \\\leq&~ \notag
            \big(1 + 2L_{1}(1+\hat{L}^{2})\big) e^{T} 
            \Delta \sum_{i=0}^{k-1} \mathbb{E}
            \big[\big|Pe_{i}\big|^{2}\big]
            +
            \frac{p_1-1}{p_1-2} e^{T} 
            \sum_{i=1}^{K} \left|\left| \mathcal{Q}_i \right|\right|^2
            \\&~+ 
            \frac{\theta^{2}+1}{\theta^{2}} e^{T} \Delta^{-1}
            \sum_{i=1}^{K} \mathbb{E}\big[\big| 
            \mathbb{E} \big(\mathcal{Q}_{i} \,|\, 
            \mathcal{F}_{t_{i-1}}\big)\big|^{2}\big].
      \end{align}
      An application of \eqref{asm:FG-FG2}, \eqref{lem:error2:final} and $\theta\Delta < \frac{1}{4L_{1}(1 + \hat{L}^{2})}$ gives that for $k = 0,1,\cdots,K$,
      \begin{align*}
            \big|Pe_k - \theta\big(f(t_{k},U_{t_k}) 
            - f(t_{k},u_{k})\big)\Delta\big|^{2}
            \geq&~
            \big|Pe_{k}\big|^{2} 
            - 
            2\theta\Delta\big\langle Pe_{k}, 
            f(t_{k},U_{t_k}) 
            - f(t_{k},u_{k})\big\rangle
            \\\geq&~
            \big|Pe_{k}\big|^{2} 
            - \theta\Delta L_{1}|e_{k}|^{2}
            \geq 
            \big(1 - 2\theta\Delta
            L_{1}(1 + \hat{L}^{2}) \big)
            \big|Pe_{k}\big|^{2}
            \\\geq&~ 
            \frac{1}{2}\big|Pe_{k}\big|^{2},
      \end{align*}
      which together with \eqref{lem:error2:leq3} yields
      \begin{align*}
            \mathbb{E}\big[\big|Pe_{k}\big|^{2}\big]
            \leq&~ \notag
            2\big(1 + 2L_{1}(1+\hat{L}^{2})\big) e^{T} 
            \Delta \sum_{i=0}^{k-1} \mathbb{E}
            \big[\big|Pe_{i}\big|^{2}\big]
            +
            2\frac{p_1-1}{p_1-2} e^{T} 
            \sum_{i=1}^{K} \left|\left| \mathcal{Q}_i \right|\right|^2
            \\&~+ 
            2\frac{\theta^{2}+1}{\theta^{2}} e^{T} 
            \Delta^{-1} \sum_{i=1}^{K} 
            \mathbb{E}\big[\big| 
            \mathbb{E} \big(\mathcal{Q}_{i} \,|\, 
            \mathcal{F}_{t_{i-1}}\big)\big|^{2}\big].
      \end{align*}
      Then the discrete Gronwall inequality indicates the desired result and completes the proof.\hfill$\square$

The lemma below concerns the coefficients of \eqref{eq:inherent_U} under Assumption \ref{asm:FG}.

\begin{lemma}\label{lem:fg2}
      Let $A_t \equiv A, t \in [0,T]$ be a constant matrix and suppose that Assumption \ref{asm:FG} holds. Then there exists a constant $C > 0$ such that for any $x,y \in \mathbb{R}^{d}$ and $s,t \in [0,T]$,
      \begin{gather*}
            \label{lem:fg2:rs1}
            \big|f(t,x)-f(s,y)\big|
            \leq
            C\Big(\big(1 + |Px| + |Py|\big)^{\gamma-1} 
            |Px-Py| + \big(1 + |Px| + |Py|\big)^{\gamma}
            |t-s|\Big),
            \\\label{lem:fg2:rs2}
            \big|g(t,x)-g(s,y)\big|^{2}
            \leq
            C\Big(\big(1 + |Px| + |Py|\big)^{\gamma-1} 
            |Px-Py|^{2} + \big(1 + |Px| 
            + |Py|\big)^{\gamma+1} |t-s|^{2}\Big).
      \end{gather*}
\end{lemma}

\textbf{Proof.}\quad
      For any $x,y \in \mathbb{R}^{d}$ and $s,t \in [0,T]$,
      \eqref{asm:FG-absF}, \eqref{lem-v-rs1} and \eqref{lem-v-rs3} imply
      \begin{align*}
            &~\big|f(t,x) - f(s,y)\big|
            =
            \big|A^-F(t,x+\hat{V}(t,x))
            - A^-F(s,y+\hat{V}(s,y))\big|
            \\\leq&~
            \big|A^-\big|C\Big( \big(1 
            + \big|x+\hat{V}(t,x)\big| 
            + \big|y+\hat{V}(s,y)\big|\big)^{\gamma-1} 
            \big|\big(x+\hat{V}(t,x)\big)
            \\&~- \big(y+\hat{V}(s,y)\big)\big|
            +
            \big(1 + \big|x+\hat{V}(t,x)\big|
            + \big|y+\hat{V}(s,y)\big|\big)^{\gamma}
            \big|t-s\big|\Big) 
            \\\leq&~
            C\Big( \big(1 + \big|Px\big| 
            + \big|Py\big|\big)^{\gamma-1} 
            \big(\big|Px-Py\big| + \big|t-s\big|\big)
            +
            \big(1 + \big|Px\big| + \big|Py\big|\big)^{\gamma} 
            \big|t-s\big|\Big)
            \\\leq&~
            C\Big(\big(1 + \big|Px\big| 
            + \big|Py\big|\big)^{\gamma-1}\big|Px-Py\big|
            +
            \big(1 + \big|Px\big| 
            + \big|Py\big|\big)^{\gamma}\big|t-s\big|\Big).
      \end{align*}
      Similarly, using \eqref{asm:FG-absG}, \eqref{lem-v-rs1} and \eqref{lem-v-rs3} again leads to 
      \begin{align*}
            &~\big|g(t,x) - g(s,y)\big|^{2}
            =
            \big|A^-G(t,x+\hat{V}(t,x))
               - A^-G(s,y+\hat{V}(s,y))\big|^{2}
            \\\leq &~
            \big|A^-\big|^{2}\big|G(t,x+\hat{V}(t,x))
            - G(s,y+\hat{V}(s,y))\big|^{2}
            \\\leq &~
            \big|A^-\big|^{2}C\Big(\big(1 
            + \big|x+\hat{V}(t,x)\big| 
            + \big|y+\hat{V}(s,y)\big|\big)^{\gamma-1} 
            \big|\big(x+\hat{V}(t,x)\big)
            \\&~
            - \big(y+\hat{V}(s,y)\big)\big|^{2}
            + \big(1 + \big|x+\hat{V}(t,x)\big|
            + \big|y+\hat{V}(s,y)\big|\big)^{\gamma+1} 
            \big|t-s\big|^{2}\Big) 
            \\\leq&~
            C\Big(\big(1 + \big|Px\big|
            + \big|Py\big|\big)^{\gamma-1} 
            \big(|Px-Py|^{2} + |t-s|^{2}\big)
            + 
            \big(1 + \big|Px\big| + \big|Py\big|)^{\gamma+1} 
            \big|t-s\big|^{2}\Big)
            \\\leq&~
            C\Big(\big(1 + \big|Px\big| 
            + \big|Py\big|\big)^{\gamma-1}\big|Px-Py\big|^{2}
            + \big(1 + \big|Px\big| + \big|Py\big|
            \big)^{\gamma+1} \big|t-s\big|^{2}\Big).
      \end{align*}
      Thus we complete the proof.\hfill$\square$

We are in a position to develop another main result on mean square convergence rate under less restrictive conditions.

\begin{theorem}\label{th:conv2A}
      Let $A_t \equiv A, t \in [0,T]$ be a constant matrix and suppose that Assumption \ref{asm:FG} holds. If $\theta \in [1/2,1]$ and $\Delta < \frac{1}{2L_{1}\theta(1 + \hat{L}^{2})}$, then there exists a constant $C > 0$, independent of $\Delta$, such that
      \begin{equation*}
            \max_{0\leq k \leq K}\mathbb{E}
            \big[\big|X_{t_k} - x_k\big|^{2}\big]
            \leq 
            C\Delta.
      \end{equation*}
\end{theorem}

\textbf{Proof.}\quad
      Noting that $U_{t} = PX_{t}$ and $PU_{t} = P^{2}X_{t} = PX_{t} = U_{t}$ for any $t \in [0,T]$, Lemma \ref{lem:fg2} and the H\"older inequlity shows that for any $s \in [t_{k-1}, t_k)$,
      \begin{align*}
            &~\mathbb{E}\big[\big|f(s,U_s) 
            - f(t_{k},U_{t_{k}})\big|^{2}\big]
            \\\leq&~
            C\Big(\mathbb{E}\big[\big(1 + \big|U_s\big|
            + \big|U_{t_{k}}\big|\big)^{2\gamma-2} 
            \big|U_s-U_{t_{k}}\big|^{2}\big]
            +
            \mathbb{E}\big[\big(1 + \big|U_s\big|
            + \big|U_{t_{k}}\big|\big)^{2\gamma}\big]
            (t_k-s)^{2}\Big)
            \\\leq&~
            C\Big(\mathbb{E}\big[\big(1 + \big|U_s\big|
            + \big|U_{t_{k}}\big|\big)^{4\gamma-2}
            \Big)^{\frac{\gamma-1}{2\gamma-1}} 
            \Big(\mathbb{E}\big[\big|U_s
            - U_{t_{k}}\big|^{\frac{4\gamma-2}{\gamma}}\big]
            \Big)^{\frac{\gamma}{2\gamma-1}} 
            \\&~+
            C\Delta^{2} \mathbb{E}\big[\big(1 + \big|U_s\big|
            + \big|U_{t_{k}}\big|\big)^{2\gamma}\big]
            \\\leq&~
            C\Delta,
      \end{align*}
      where $2\gamma \leq 4\gamma-2 < p_{1}$ and Lemmas \ref{lem:bound}, \ref{lem:UX} have been used. Similarly, one can prove that for any $s \in [t_{k-1},t_{k})$,
      \begin{gather*}
            \label{th:conv2:f2}
            \mathbb{E}\big[\big|f(s,U_s) 
            - f(t_{k-1},U_{t_{k-1}})\big|^{2}\big]
            \leq 
            C\Delta,
            \quad
            \mathbb{E}\big[\big|g(s,U_s) 
            - g(t_{k-1},U_{t_{k-1}})\big|^{2}\big]
            \leq 
            C\Delta.
      \end{gather*}
      Together with the techniques used in \eqref{lem:conv:r1} and \eqref{lem:conv:r2}, we analogously obtain that for any $k = 1,2, \cdots, K$,
      \begin{equation*}
             \mathbb{E}\big[\big|\mathcal{Q}_{k}\big|^{2}\big]
             \leq
             C\Delta^{2},
             \quad
             \mathbb{E}\big[\big|\mathbb{E}\big( 
             \mathcal{Q}_{k} \,|\, \mathcal{F}_{t_{k-1}} 
             \big) \big|^{2}\big]
             \leq
             C\Delta^{3}.
      \end{equation*}
      Combining Lemma \ref{lem:error2} yields the desired result and thus ends the proof.\hfill$\square$

\section{Numerical experiments}\label{sec:experiments}
This section is devoted to carrying out numerical experiments that illustrate the theoretical results obtained earlier. To begin with, we emphasize that in all subsequent simulations, the nonlinear equations arising from the implementation of the implicit method at each time step are solved using Newton–Raphson iterations with a tolerance of $10^{-5}$. To examine the mean square convergence rate of each STM, the unavailable exact solution of a given SDAEs is approximated by a reference solution, which is generated by \eqref{eq:ABEM} with a sufficiently fine stepsize $\Delta = 2^{-13}$. The other numerical approximations are computed by applying \eqref{eq:ABEM} with six coarser stepsizes $\Delta = 2^{-i}, i = 6, 7, 8, 9, 10, 11$. Furthermore, all expectations are estimated through Monte Carlo simulations based on $1000$ independent Brownian sample paths.

\begin{example}
      \rm
      Consider the following $\R^2$-valued index 1 SDAEs
      \begin{equation}\label{eq:exam2}
            A_t \diff X_t 
            = 
            F(t,X_t) \text{d}t + G(t,X_t) \diff W_t,
            \quad t \in [0,T],
      \end{equation}
      where $\{W_{t}\}_{t \in [0,T]}$ is a standard $\mathbb{R}^{2}$-valued Brownian motion. Let $A_t = \begin{pmatrix} 0 & 0 \\ -\frac{1}{\sqrt{2}}(t^2+1) & \frac{1}{\sqrt{2}}(t^2+1) \end{pmatrix}$ be a singular matrix with the singular value decomposition
      \begin{equation*}
            A_{t} 
            =
            \begin{pmatrix} 0 & 1 \\ 1 & 0 \end{pmatrix}
            \begin{pmatrix} t^2+1 & 0 \\ 0 & 0 \end{pmatrix}
            \begin{pmatrix} -\frac{1}{\sqrt{2}} & \frac{1}{\sqrt{2}} 
            \\ \frac{1}{\sqrt{2}} & \frac{1}{\sqrt{2}} \end{pmatrix}
            = M \Sigma_t N .
      \end{equation*}
      It follows that 
      \begin{equation*}
            A_t^- 
            = 
            \begin{pmatrix} 0 & -\tfrac{1}{\sqrt{2}(t^2+1)} 
            \\ 0 & \tfrac{1}{\sqrt{2}(t^2+1)} \end{pmatrix},
            \quad
            P = \begin{pmatrix} \frac{1}{2} & -\frac{1}{2} 
            \\ -\frac{1}{2} & \frac{1}{2} \end{pmatrix},
            \quad
            R = \begin{pmatrix} 1 & 0 
            \\ 0 & 0 \end{pmatrix}.
      \end{equation*}
      For any $x \in \mathbb{R}^{2}$ and $t \in [0,T]$, let
      \begin{equation*}
            F(t,x) 
            = 
            \begin{pmatrix} x_1 + x_2 + \sin(t) 
            \\ (x_1-x_2)^3-a(x_1-x_2)+1 \end{pmatrix}, 
            \quad
            G(t,x)
            =
            \begin{pmatrix} 0 & 0 
            \\ b(x_1+x_2+1) 
            & b(x_1-x_2)^2 +2b \end{pmatrix}.
      \end{equation*}
      Thus, requiring $(p_{1}-1)b^{2} \leq \tfrac{3}{4}$ and and using $m^2+mn+n^2 \geq \tfrac{3}{4}(m+n)^2,m,n \in \mathbb{R}$ shows that for any $x,y \in \mathbb{R}^{2}$ and $t \in [0,T]$,
      \begin{align}\label{eq:exam2:PxAF}
            &~2\big\langle Px - Py, A_t^- F(t,x) 
            - A_t^- F(t,y) \big\rangle 
            + 
            (p_1-1)\big|A_t^-G(t,x) 
            - A_t^-G(t,y)\big|^{2} \notag
        \\=&~
            -\frac{\sqrt{2}}{t^2+1} 
            \big(x_1-x_2-y_1+y_2\big)^{2} 
            \Big(\big(x_1-x_2\big)^{2} 
            + \big(x_1-x_2\big)\big(y_1-y_2\big) 
            + \big(y_1-y_2\big)^{2} -a \Big) \notag
            \\&~
            + \frac{(p_1-1)b^2}{t^2+1} 
            \Big(\big(x_1+x_2-y_1-y_2\big)^{2} 
            + \big(x_1-x_2-y_1+y_2\big)^{2}
            \big(x_1-x_2+y_1-y_2\big)^{2}\Big) \notag
            \\\leq&~
            -\frac{\sqrt{2}}{t^2+1}
            \big(x_1-x_2-y_1+y_2\big)^{2} 
            \Big( \big(x_1-x_2\big)^2 
            + \big(x_1-x_2\big)\big(y_1-y_2\big) 
            + \big(y_1-y_2\big)^{2}   \notag
            \\&~- 
            (p_1-1)b^{2}\big(x_1-x_2
            +y_1-y_2\big)^{2}\Big)
            +
            \frac{\sqrt{2}}{t^2+1}
            \Big(a(x_1-x_2-y_1+y_2)^2 \notag
            \\&~ 
            + (p_1-1)b^{2}\big(x_1
            + x_2-y_1-y_2\big)^{2}\Big) \notag
            \\\leq&~ 
            -\frac{\sqrt{2}}{t^2+1}
            \big(x_1-x_2-y_1+y_2\big)^{2} 
            \big(\tfrac{3}{4} -(p_1-1)b^{2}\big)
            \big(x_1-x_2+y_1-y_2\big)^{2} \notag
            \\&~ 
            +\frac{2\sqrt{2}}{t^2+1}
            \Big( a(x_1-y_1)^2+a(x_2-y_2)^2
            +(p_1-1)b^{2}\big(x_1-y_1\big)^{2}
            +(p_1-1)b^{2}\big(x_2-y_2\big)^{2}\Big) \notag
            \\\leq&~
            \frac{2\sqrt{2}}{t^2+1} \big(a + (p_1-1)b^{2} \big) 
            \big( (x_1-y_1)^2 + (x_2-y_2)^2 \big) \notag 
            \\\leq&~
            2\big(a + (p_{1}-1)b^{2}\big)|x-y|^2, \notag
      \end{align}
      and
      \begin{align}
            \big|F(t,x) - F(s,y)\big|
            \leq&~
            \big|x_1+x_2+\sin(t)
            - y_1-y_2-\sin(s)\big| \notag
            \\&~+
            \big|(x_1-x_2)^3-a(x_1-x_2)
            -(y_1-y_2)^3+a(y_1-y_2)\big| \notag
            \\\leq &~
            \big|(x_1-x_2)^3 - (y_1-y_2)^3\big| 
            + C|x-y| + |t-s| \notag
            \\\leq&~
            C\big(1+|x|+|y|\big)^{2}|x-y| 
            + |t-s|, \notag
      \end{align}
      as well as 
      \begin{align*}
            \big|G(t,x) - G(s,y)\big|^{2}
            =&~
            b^{2}(x_1+x_2-y_1-y_2)^2 
            +
            b^{2}((x_1-x_2)^2-(y_1-y_2)^2)^2
            \\\leq&~
            C (1+|x|+|y|)^{2} |x-y|^2.
      \end{align*}
      These inequalities imply (A1) in Assumption \ref{asm:FG}. In view of
      \begin{align*}
            J(t,x) = A_t + RF'_x(t,x)
            = 
            \begin{pmatrix} 1 & 1 
            \\ -\tfrac{t^2+1}{\sqrt{2}}
            &   \tfrac{t^2+1}{\sqrt{2}} \end{pmatrix},
      \end{align*}
      we get
      \begin{align*}
            \big|J(t,x)^{-1}\big|
            = 
            \left|\begin{pmatrix} \tfrac{1}{2} 
            & -\tfrac{\sqrt{2}}{2(t^2+1)} \\ 
            \tfrac{1}{2} & 
            \tfrac{\sqrt{2}}{2(t^2+1)} 
            \end{pmatrix}\right|
            =
            \bigg(\frac{1}{2} 
            + \frac{1}{(t^2+1)^{2}}\bigg)^{\frac{1}{2}}
            \leq 
            \sqrt{\frac{3}{2}}.
      \end{align*}
      Thus (A2) in Assumption \ref{asm:FG} holds. Let the initial value $X_0=(1, -1)$, then $RF(0,X_0) = 1 - 1 + \sin(0) = 0$, which implies (A3) in Assumption \ref{asm:FG}. Noting that Assumption \ref{asm:At} holds with $r=1$, applying Corollary
      \ref{co:conv} implies that each STM \eqref{eq:ABEM} for \eqref{eq:exam2} converges strongly with order $\frac{1}{2}$.

\begin{figure}[!htbp]
\begin{center}
      \subfigure[$\theta = 0.5$]{
      \includegraphics[width=0.32\textwidth]{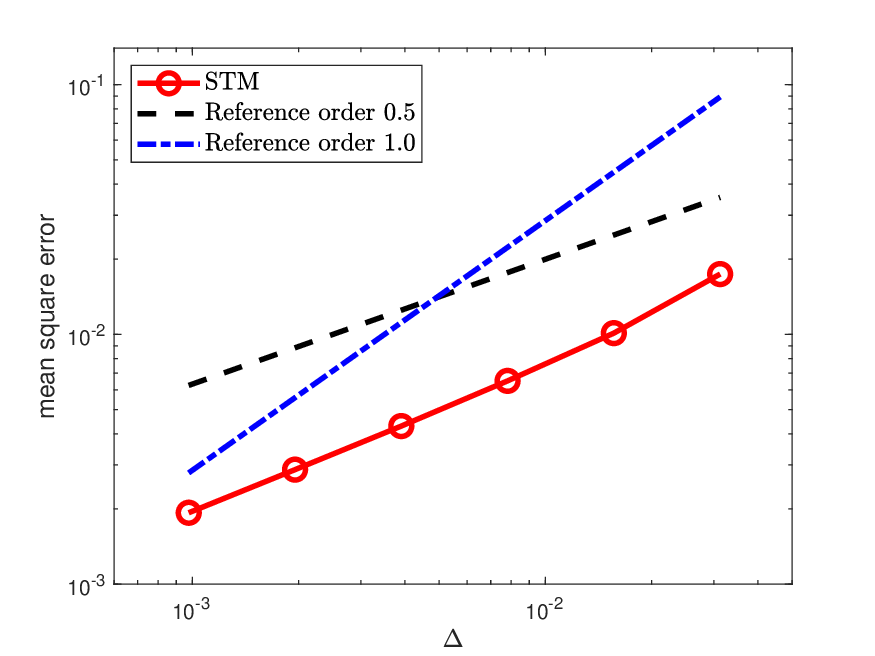}}
      \subfigure[$\theta = 0.75$]{
      \includegraphics[width=0.32\textwidth]{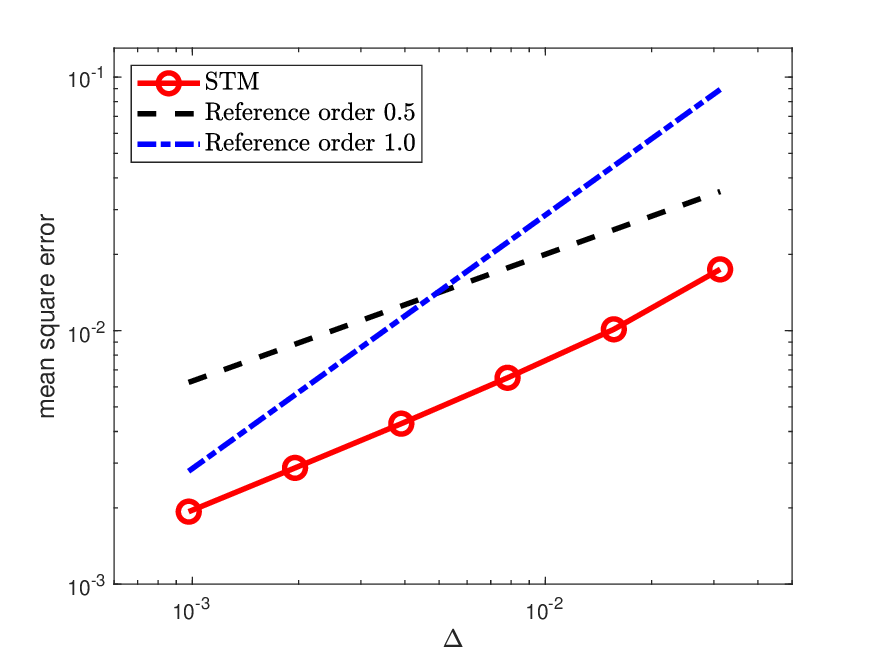}}
      \subfigure[$\theta = 1.0$]{
      \includegraphics[width=0.32\textwidth]{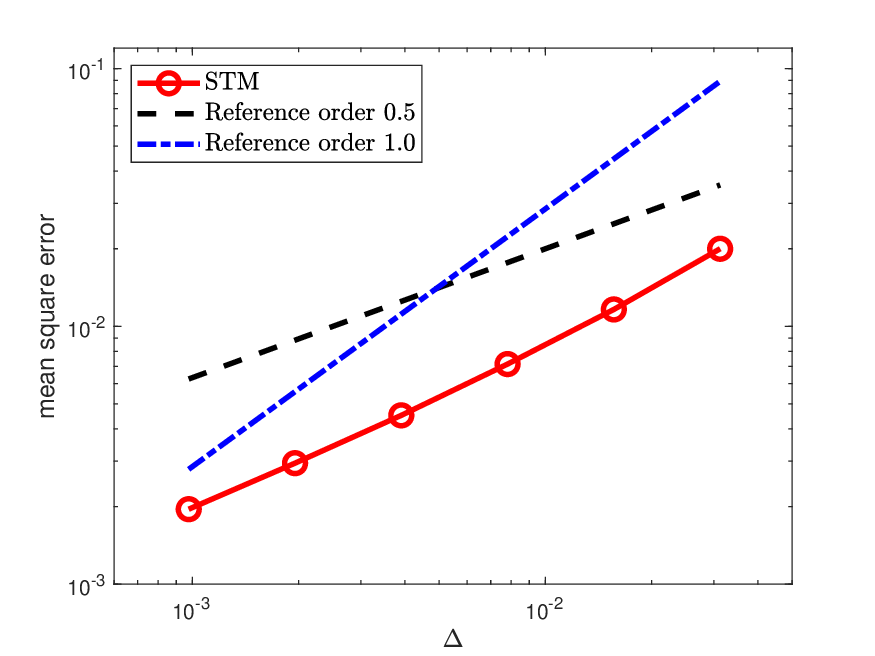}}
      \centering\caption{Mean square convergence rates of STMs applied to \eqref{eq:exam2}}
      \label{fig:2Dorder}
\end{center}
\end{figure}

      In view of $\gamma = 3$, $p_{1} > 4\gamma-2$ and $(p_{1}-1)b^{2} \leq \tfrac{3}{4}$, we take $a = 1$ and $b = \frac{1}{5}$. 
      From Figure \ref{fig:2Dorder}, we see that the root mean square error lines and the reference line appear to parallel to each other for different choices of $\theta$, showing that the mean square convergence rate of each STM is order $\frac{1}{2}$. Moreover, a least square fit indicates that the slope of the line for STM is $0.6264$ for $\theta = 0.5$, $0.6264$ for $\theta = 0.75$ and $0.6681$ for $\theta = 1.0$, respectively, which identifies Corollary \ref{co:conv} numerically.
\end{example}

\begin{example}
      \rm
      Consider the following $\R^3$-valued index 1 SDAEs as follows
      \begin{equation}\label{eq:exam1}
            A_t \diff X_t 
            = F(t,X_t) \text{d}t + G(t,X_t) \diff W_t,
            \quad t \in [0,T]. 
      \end{equation}
      where $\{W_{t}\}_{t \in [0,T]}$ is a standard $\mathbb{R}^{3}$-valued Brownian motion. Let $A_t = \mathrm{diag}\big(\frac{1}{2(t^2+1)}, 10, 0\big)$ be a singular matrix, then Assumption \ref{asm:At} holds with $r = 2$, $\underline{\sigma} = \frac{1}{2(T^2+1)}$ and $\overline{\sigma} = 10$. It can be calculated that $A_t^- = \mathrm{diag}\big(2(t^2+1), \tfrac{1}{10}, 0\big)$, $P = \mathrm{diag}(1, 1, 0)$, and $R = \mathrm{diag}(0, 0, 1)$. For any $x \in \mathbb{R}^{3}$ and $t \in [0,T]$, let
      \begin{equation*}
            F(t,x) 
            = 
            \big(-x_1^3, x_3, x_2t+x_3\big)^{\top}, 
            \quad
            G(t,x)= 
            \mathrm{diag}(\sin(t), c x_1^2, 0).
      \end{equation*}
      By requiring $(p_{1}-1)c^{2} \leq 300$ and using $m^2 + mn + n^2 \geq \tfrac{3}{4}(m+n)^2, m,n \in \mathbb{R}$, we see that for $x,y \in \R^3$ and $t \in [0,T]$,
      \begin{align}
            &~2\big\langle Px -Py, A^-_t F(t,x)   \notag
            - A^-_t F(t,y) \big\rangle 
            + 
            (p_1-1)\big|A^-_t G(t,x) - A^-_t G(t,y)\big|^{2}
            \\=&~ \notag
            -4\big(t^{2}+1\big)\big(x_1-y_1\big)
            \big(x_1^3-y_1^3\big) 
            +
            \frac{1}{5}\big(x_2-y_2\big)\big(x_3-y_3\big) 
            +
            (p_1-1)\bigg(\frac{c}{10}
            \big(x_1^2-y_1^2\big)\bigg)^{2}
            \\\leq&~ \notag
            -4\big(x_1-y_1\big)^{2}
            \big(x_1^2+x_1y_1+y_1^2\big) 
            +
            \frac{p_1-1}{100}c^{2}
            \big(x_1-y_1\big)^{2}\big(x_1+y_1\big)^{2} 
            +
            \frac{(x_2-y_2)^{2} +(x_3-y_3)^{2}}{10} 
            \\\leq&~ \notag
            -\big(x_1-y_1\big)^{2}\bigg(
            4\big(x_1^2+x_1y_1+y_1^2\big) 
            - \frac{p_1-1}{100}c^{2}
              \big(x_1+y_1\big)^{2} \bigg) 
            + 
            \frac{1}{10}|x-y|^{2} 
            \\\leq&~ \notag 
            -\bigg(3 - \frac{p_1-1}{100}c^{2} \bigg)
            \big(x_1-y_1\big)^{2}\big(x_1+y_1\big)^2 
            + 
            \frac{1}{10}|x-y|^{2}
            \\\leq&~ \notag
            \frac{1}{10}|x-y|^{2},
      \end{align} 
      and
      \begin{align}
            \big|F(t,x) - F(s,y)\big|
            =&~ \notag
            \big((x_1^3-y_1^3)^2 +(x_3-y_3)^2 
            + (x_2t+x_3 -y_2s-y_3)^2\big)^{\frac{1}{2}}
            \\\leq&~ \notag
            |x_1^3-y_1^3| + 2|x_3 -y_3| 
            + |x_2t -y_2t| + |y_2t -y_2s|
            \\\leq&~ \notag
            C\big(1 + |x| + |y|\big)^{2}|x-y| 
            + |y||t-s|
            \\\leq&~\notag
            C\Big(\big(1 + |x| + |y|\big)^{2}|x-y| 
            + \big(1 + |x| + |y|\big)^{3}|t-s|\Big),
      \end{align}
      as well as
      \begin{align}\label{eq:exam1:G-G}
            \big|G(t,x) - G(s,y)\big|^{2}
            =&~\notag
            \big(\sin(t) - \sin(s)\big)^{2} 
            +
            b^{2}\big(x_1^2-y_1^2\big)^{2}
            \\\leq&~ \notag
            C\big(1+|x|+|y|\big)^{2}|x-y|^{2} + |t-s|^{2}
            \\\leq&~ \notag
            C\Big(\big(1+|x|+|y|\big)^{2}|x-y|^{2} 
            +  \big(1 + |x| + |y|\big)^{4}|t-s|^{2}\Big).
      \end{align}
      Then (A1) in Assumption \ref{asm:FG} holds with $L_{1} = \frac{1}{10}$, $\gamma = 3$ and $p_{1} > 10$. Noting that 
      \begin{align*}
            J(t,x) = A_t + RF'_x(t,x)
            = \begin{pmatrix} \tfrac{1}{2(t^2+1)} & 0 & 0 
            \\ 0 & 10 & 0 \\ 0 & t & 1 \end{pmatrix},
            \quad x \in \mathbb{R}^{d}, t \in [0,T],
      \end{align*}
      we have
      \begin{align*}
            J(t,x)^{-1}
            = \begin{pmatrix} 2(t^2+1) & 0 & 0 
            \\ 0 & \tfrac{1}{10} & 0 
            \\ 0 & -\tfrac{t}{10} & 1 \end{pmatrix},
      \end{align*}
      and consequently
      \begin{equation*}
            \big|J(t,x)^{-1}\big| 
            = 
            \big( 4(t^2+1)^{2} + 0.01 
            + 0.01t^{2} + 1\big)^{\frac{1}{2}} 
            \leq 
            2(T^{2}+1) + 2(T+1).
      \end{equation*}
      Thus (A2) in Assumption \ref{asm:FG} holds. Letting the initial value $X_0 = (1, -1, 0)$ yields $RF(0,X_0) = 0$, which implies (A3) in Assumption \ref{asm:FG}. Therefore, we conclude that Assumption \ref{asm:FG} holds. Besides, imposing the condition $p_{2}c^{2} \leq \frac{3}{4(T^{2}+1)}$ enables us to calculate that for any $x,y \in \mathbb{R}^{d}$ and $t \in [0,T]$, 
      \begin{align}\label{eq:exam1:AxF}
            &~2\big\langle A_tx - A_ty, 
            F(t,x) - F(t,y) \big\rangle 
            + 
            p_{2}\big|G(t,x) - G(t,y)\big|^{2} \notag
            \\=&~
            -\frac{1}{t^2+1}\big(x_1-y_1\big)
            \big(x_1^3-y_1^3\big) 
            +
            20\big(x_2-y_2\big)\big(x_3-y_3\big) 
            +
            p_{2} c^{2}\big(x_1^2-y_1^2\big)^{2} \notag
            \\\leq&~
            -\big(x_1-y_1\big)^{2} 
            \bigg(\frac{x_1^2+x_1y_1+y_1^2}{T^2+1}
            - p_{2}c^{2}\big(x_1+y_1\big)^{2}\bigg) 
            + 
            10|x_2-y_2|^{2} +10|x_3-y_3|^{2} \notag
            \\\leq&~
            -\bigg(\frac{3}{4(T^2+1)}-p_2c^2\bigg)
            \big(x_1-y_1\big)^{2}\big(x_1+y_1\big)^2 
            + 
            10\big|x-y\big|^{2} \notag
            \\\leq&~ 
            10\big|x-y\big|^{2}, \notag
      \end{align}
      which means that Assumption \ref{asm:FGAt} holds with $L_{2} = 10$. As a result, Theorem \ref{thm:convAt} holds for \eqref{eq:exam1}. Concerning $\gamma = 3$, $p_{1} > 4\gamma-2$, $(p_{1}-1)c^{2} \leq 300$ and $p_{2}c^{2} \leq \frac{3}{4(T^{2}+1)}$, we take $T = 1$, $p_{1} = 21$, $p_{2} = 25$ and $c = \frac{1}{10}$. In Figure \ref{fig:2Dorder}, the root mean square error lines and the reference line with order $1.0$ are nearly parallel, which means that the mean square convergence rate of each STM is $1.0$. It is beyond the theoretical order of convergence for this example.  

\begin{figure}[!htbp]
\begin{center}
      \subfigure[$\theta = 0.5$]{
      \includegraphics[width=0.32\textwidth]{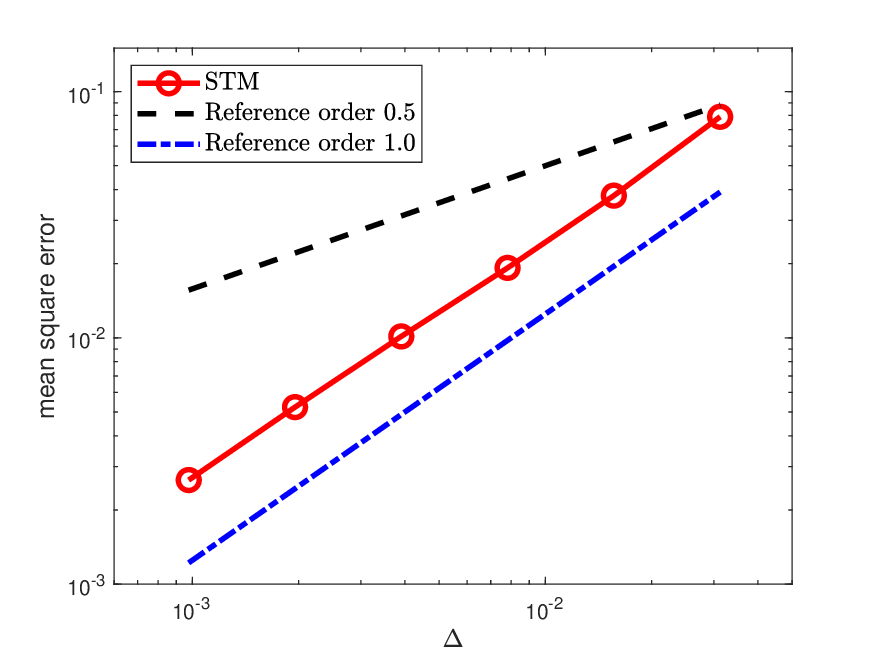}}
      \subfigure[$\theta = 0.75$]{
      \includegraphics[width=0.32\textwidth]{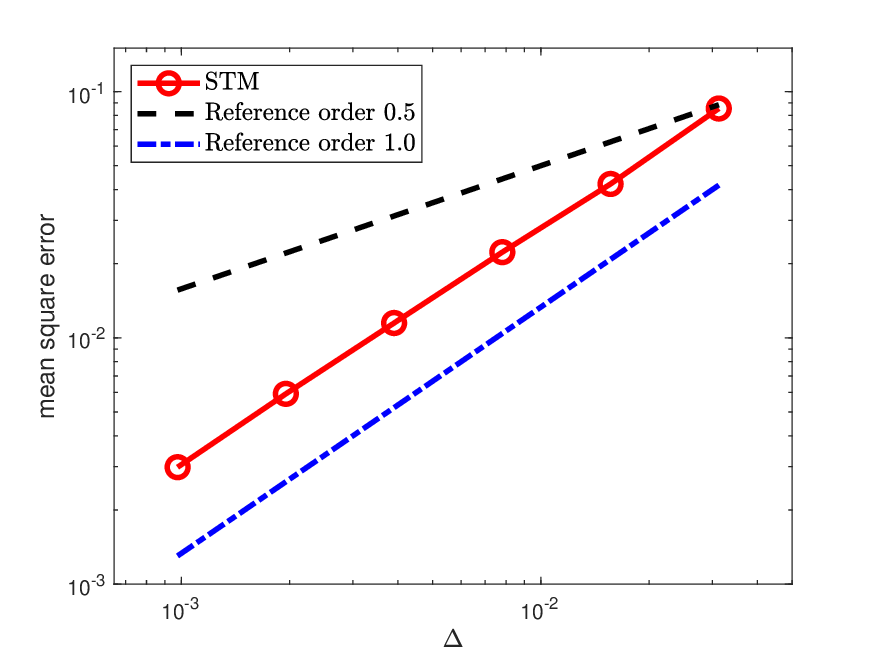}}
      \subfigure[$\theta = 1.0$]{
      \includegraphics[width=0.32\textwidth]{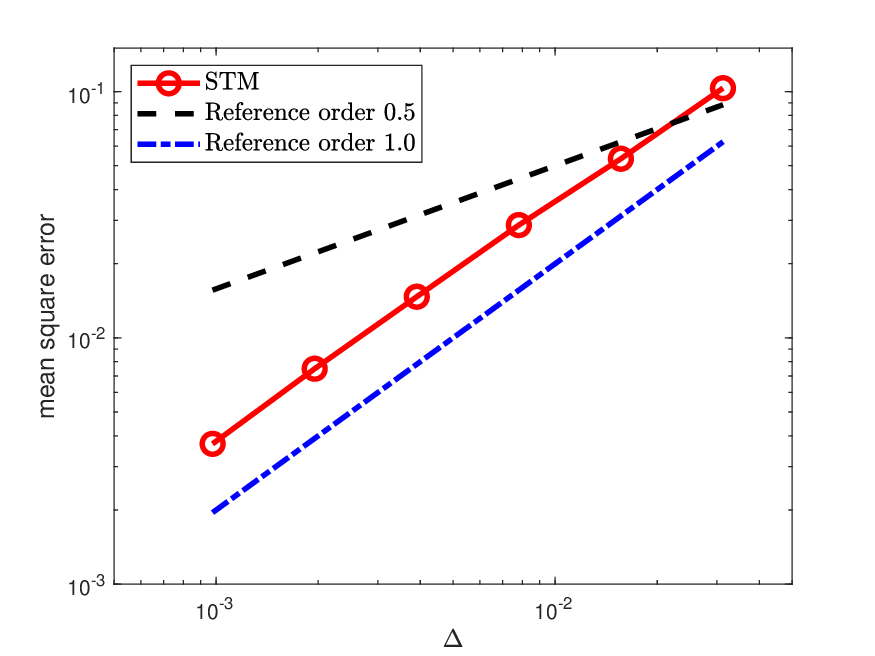}}
      \centering\caption{Mean square convergence rates of STMs applied to \eqref{eq:exam1}}
      \label{fig:3Dorder}
\end{center}
\end{figure}
     
\end{example}

\section{Conclusion}
This work investigates the strong convergence rate of STMs for index 1 SDAEs with non-constant singular matrices under non-global Lipschitz conditions. We establish that each STM with $\theta \in [\frac{1}{2},1]$ attains a mean square convergence rate of order $\frac{1}{2}$, and these theoretical results are further validated by numerical experiments. Building upon the convergence analysis framework developed here, a natural direction for future research is the extension to Milstein-type methods. Moreover, it is also of great interest to study numerical approximations for SDAEs of index 2 \cite{the2024stochastic}, which arise frequently in physical systems with strong constraints, such as multi-body dynamics and degradation models of circuit systems. However, the higher index generally leads to weaker regularity of solutions, which prevents a direct extension of the present analysis. Addressing these challenges will require the development of new techniques and arguments.

\bibliographystyle{plain}	

\end{document}